\documentclass[reqno,11pt]{amsart}
\usepackage{amsmath,amssymb,amsthm,amsfonts,mathrsfs,times,color}
\usepackage{amsmath,amssymb,amsthm,amsfonts,mathrsfs,latexsym,
times,color,epsfig,bm}
\usepackage{ulem}
\usepackage{amsmath,amsfonts,mathtools}
\usepackage{amsthm}
\usepackage{graphicx}
\usepackage{todonotes}
\usepackage{color}
\usepackage{subcaption}
\usepackage{comment}

\usepackage{multicol}
\newcommand{\padi}[1]{\textcolor{black}{#1}}

\newtheorem{thm}{Theorem}[section]

\newtheorem{lemma}[thm]{Lemma}
\newtheorem{remark}{Remark}[section]
\newtheorem{theorem}[thm]{Theorem}

\numberwithin{equation}{section}

\newcommand{\tu}{{\tilde u}}

\newcommand{\tv}{{\tilde v}}

\newcommand{\dx}{{\mathrm{d}x}}

\usepackage{epsfig}
\usepackage{graphicx}
\usepackage{colordvi}
\usepackage{graphics}
\usepackage{color}
\usepackage{multirow}

\makeatletter
\@namedef{subjclassname@2020}{%
  \textup{2020} Mathematics Subject Classification}
\makeatother

\begin{document}

\title[On a Logarithmically Sensitive Chemotaxis Model]{Global stability of a logarithmically sensitive chemotaxis model under time-dependent  boundary conditions}

\author[P. Fuster Aguilera]{Padi Fuster Aguilera}
\address[P. Fuster Aguilera]{Department of Mathematics, University of Colorado Boulder, Boulder, CO 80309, USA}
\email{padi.fuster@colorado.edu}

\author[K. Zhao]{Kun Zhao}
\address[K. Zhao]{Department of Mathematics, Tulane University, New Orleans, LA 70118, USA}
\email{kzhao@tulane.edu}


\begin{abstract}
This paper studies the dynamical behavior of classical solutions to a hyperbolic system of balance laws, derived from a chemotaxis model with logarithmic sensitivity, subject to time-dependent boundary conditions. It is shown that under suitable assumptions on the boundary data, solutions starting in $H^2$-space exist globally in time and the differences between the solutions and their corresponding boundary data converge to zero, as time goes to infinity. There is no smallness restriction on the magnitude of initial perturbations. Moreover, numerical simulations show that the assumptions on the boundary data are necessary for the above mentioned results.
\end{abstract}
\maketitle

\section{Introduction}
The movement of an organism or entity in response to a chemical stimulus in the local environment is termed as chemotaxis. It is the underlying mechanism of many biological processes in modern cell biology, biochemistry and clinical pathology, such as tumor angiogenesis \cite{angiogenesis}, blood vessel formation \cite{blood-vessel}, slime mould formation \cite{slime}, bacterial foraging, immune response, embryonic development and tissue homeostasis \cite{Murray,TLM}, fish pigmentation patterning \cite{fish}, primitive streak formation \cite{primitive}, wound healing \cite{wound}, just to mention a few. Early scholarly descriptions of chemotaxis were introduced in the late 19th century by Engelmann (1881),  Pfeffer (1884), Metchnikoff (1882-1886) and Jennings (1906). Later important contributions include, but are not limited to, quality control of chemotaxis assays ( Harris 1953) and study of intracellular signal transduction of bacteria (Adler 1966).

Mathematical modeling of chemotaxis using continuum partial differential equations was initiated in the 1950s by Patlak \cite{Patlak} from a probabilistic perspective, and regained popularity in the 1970s through the pioneering work of Keller and Segel \cite{KS1,KS2,KS3} using a phenomenological approach. The original Keller-Segel models contain two prototypes according to the chemotactic sensitivity function - linear sensitivity and logarithmic sensitivity. The linear sensitivity was employed in the modeling of self-aggregation of {\it Dictyostelium discoideum} in response to cyclic adenosine monophosphate (cAMP) \cite{KS2}, while the logarithmic sensitivity appeared in \cite{KS3} to interpret J. Adler's experimental result \cite{Adler} on the formation of traveling bands in nutrient-enticed {\it E. Coli} population.

In general form, the Keller-Segel model reads as
\begin{subequations}\label{KS}
 \begin{alignat}{2}
 u_{t} &=D u_{xx} - \chi\big[u\Phi(c)_x\big]_x,\\
 c_{t} &=\varepsilon c_{xx} - \mu u c^m - \sigma c,
 \end{alignat}
\end{subequations}
where the unknown functions $u(x,t)$ and  $c(x,t)$ denote the density of the organic population and concentration of the chemical signal at position $x$ at time $t$, respectively. The parameter $D>0$ denotes the diffusion coefficient of the organic population density; $\chi\neq0$ is the coefficient of chemotactic sensitivity, the sign of $\chi$ dictates whether the chemotaxis is attractive ($\chi>0$) or repulsive ($\chi<0$), with $|\chi|$ measuring the strength of the chemotactic response; $\varepsilon\ge0$ is the diffusion coefficient of the chemical signal; $\mu\neq0$ is the coefficient of the density-dependent production/degradation rate of the chemical signal;  $\sigma\ge0$ is the natural degradation rate of the chemical signal; and $m\ge0$ characterizes the mode of temporal growth/decay of the chemical signal. Moreover, $\Phi(c)$ denotes the chemotactic sensitivity with either $\Phi(c)=c$ (i.e., linear) or $\Phi(c)=\log(c)$ (i.e., logarithmic). The chemotactic sensitivity underlines the main character of such a model, whose spatial derivative depicts the mechanistic feature of chemotactic movement -- advection of the organic population induced by the spatial gradient of the chemical signal in the local environment. The pioneering work of Keller and Segel inspired many of the modern studies in chemotaxis research, with the frequency of publication continuing to increase.

In the direction of mathematical analysis of the two prototypes of the Keller-Segel model, a search in the database shows that, comparing with the extensive results on the model with linear sensitivity (see e.g. \cite{survey1,survey2,survey3}), much less is known for the model with logarithmic sensitivity, due to its possible singular nature. However, the logarithmic sensitivity obeys the Weber-Fechner's law, which is a fundamental principle in psychophysics and has prominent applications in biology (c.f. \cite{Weber-Fechner1,Weber-Fechner2,Weber-Fechner3,Weber-Fechner4,KS3}). The current paper concentrates upon the model with logarithmic sensitivity:
\begin{subequations}\label{OS}
 \begin{alignat}{2}
 u_{t} &=D u_{xx} - \chi\big[u(\log c)_x\big]_x,\\
 c_{t} &=\varepsilon c_{xx} - \mu u c^m - \sigma c.
 \end{alignat}
\end{subequations}
This model has been utilized in a variety of contexts to explain the underlying mechanisms of different chemotactic processes, depending on the sign of $\chi$ and $\mu$. For example, the original Keller-Segel model ($\chi>0$, $\mu>0$, $0\le m<1$, $\sigma=0$) was proposed in \cite{KS3} to describe the propagation of wave bands observed in Adler's experiment \cite{Adler}. The same model with $m=1$ was employed by Levine {\it et al} \cite{LSN} to interpret the dynamical interactions between vascular endothelial cells (VECs) and signaling molecules vascular endothelial growth factor (VEGF) in the onset of tumor angiogenesis. On the other hand, when $\chi<0$, $\mu<0$, $m=1$ and $\sigma>0$, the model was designed in \cite{LS,OS} to illustrate the chemotactic movement of reinforced random walkers (e.g.\,surface or matrix-bound adhesive molecules) that deposit non-diffusive ($\varepsilon=0$) or slowly moving ($0< \varepsilon\ll 1$) chemical signals that modify the local environment for succeeding passages.

Though the logarithmic sensitivity has been utilized to explain the mechanisms of chemotactic movement for a variety of biological phenomena, its singular nature presents a significant challenge to the qualitative analysis of the model. One supporting evidence is that except for some instability results \cite{instability1,instability2} and stability of spike-layer solutions \cite{CLW}, the stability of traveling wave solutions of \eqref{OS} with $m\neq 1$ is still an outstanding open problem in the area. On the other hand, it has been recognized that when $m=1$ the singular nature can be removed by the Cole-Hopf type transformation: $v = \frac{c_x}{c}$. Indeed, when $m=1$, by applying such a transformation, manipulating the second equation of \eqref{OS},  applying the rescalings: $\tilde{t}= |\chi\mu|D^{-1}\,t$, $\tilde{x}= \sqrt{|\chi\mu|}\,D^{-1}\,x$, $\tilde{v}= \mathrm{sign}(\chi)\sqrt{|\chi|\,|\mu|^{-1}}\ v$ and dropping the tilde, one obtains the transformed version of \eqref{OS}:
\begin{subequations}\label{TOS1}
\begin{alignat}{2}
u_t - (u v)_x &= u_{xx},\label{TOS1a}\\
v_t - \mathrm{sign}(\chi\mu) u_x &= \varepsilon{D}^{-1} v_{xx} + \varepsilon\chi^{-1} (v^2)_x,\label{TOS1b}
\end{alignat}
\end{subequations}
which is a system of dissipative conservation laws. A direct calculation shows that the eigenvalues of the Jacobian matrix associated with the flux on the left of \eqref{TOS1} are given by
\begin{align}\label{CF}
\lambda_\pm = \frac{\left(2\varepsilon\chi^{-1} -1\right)v \pm \sqrt{\left(2\varepsilon\chi^{-1}+1\right)^2v^2 + 4\,\mathrm{sign}(\chi\mu) u}}{2}.
\end{align}
This indicates that the principle part of system \eqref{TOS1} is hyperbolic in biologically relevant regimes where the cellular density $u>0$, provided $\chi\mu>0$. This enables the adaptation of fundamental analytic tools in hyperbolic balance laws, such as entropy method, to study the qualitative behavior of the model. Throughout this paper, we focus on the case of $\chi\mu>0$. It is worth mentioning that when $\chi\mu<0$, the characteristic fields may change type, which could alter the dynamics of the model drastically. This is supported by the blowup (explicit and numerical) solutions constructed in \cite{LS}.

Since the work of Levine and Sleeman \cite{LS}, mathematical study of \eqref{TOS1} has developed into an active area in applied analysis. To put things into perspective, we briefly survey the literature in connection with \eqref{TOS1} when $\chi\mu>0$. First of all, the transformed system \eqref{TOS1} satisfies the Shizuta-Kawashima condition \cite{Shizuta-Kawashima} for $\varepsilon\ge0$, which guarantees the global well-posedness and stability of classical solutions to the Cauchy problem \emph{near} constant equilibrium states (see e.g. \cite{Zeng1,Zeng2}). For large amplitude solutions, the global well-posedness and local stability results were first established in \cite{FFH,GXZZ} and \cite{ZZ}, respectively. These were upgraded in a series of recent works \cite{1d1,1d2,1d3,1d4,1d5,1d6,1d7}, where the \emph{global} stability of constant equilibrium states is demonstrated, which suggested that uniform distribution is a generic phenomenon in the process of logarithmically sensitive chemoattraction with chemical consumption or chemorepulsion with chemical production. In addition to global dynamics, the zero chemical diffusivity limit (i.e., as $\varepsilon\to0$) and instantaneous spatial analyticity of large amplitude solutions have been analyzed in \cite{DL2,DL1,1d2,1d4,1d5,1d7}. Another major direction of research is concerned with the existence and stability of traveling wave solutions of \eqref{TOS1}. In particular, the local stability of large-strength travel wave solutions has been studied in \cite{TW1,TW2,TW3,TW4,TW5,TW6,TW7,TW8,Wang-survey}. In the multi-dimensional spaces, the qualitative behaviors of the model have been investigated under various smallness assumptions on the initial data. We refer the readers to \cite{mTW1,mTW2,md8,md1,md9,md2,mdw1,md3,md4,md5,md6,md7}. Moreover, the appended version of \eqref{TOS1} with logistic growth  has been studied in \cite{FMZ,ZZ1,ZZ2,ZZ3,ZZ4} where the enhanced dissipation induced by logistic damping is explored. Furthermore, the amended model with nonlinear density-dependent chemical production/consumption rate (i.e., replacing $u_x$ by $(u^\gamma)_x$ with $\gamma>1$) has been analyzed in \cite{FXXZ,ZLMZ,ZLWZ}.

We emphasize that all of the foregoing studies of \eqref{TOS1} on a finite interval, say, $(a,b)$, subject to Dirichlet type boundary conditions require the boundary conditions to match, i.e., $u(a,t)=u(b,t)$, $v(a,t)=v(b,t)$, whether constant-valued \cite{1d2} or evolving in time \cite{1d5}. Although those results are mathematically meaningful, they are less realistic from the point of view of physical/biological applications, since in real-world situations the values of the quantity under consideration at the endpoints may be different time by time. Rigorous mathematical study of the model subject to \emph{unmatched} boundary conditions thus becomes relevant. This is the major motivation of the current paper. In addition, the model with $\varepsilon=0$ is a conceptual idealization that was designed to simplify the underlying mathematical analysis \cite{LS}. Though it has been analyzed under various types of boundary conditions, the large-time behavior of classical solutions subject to time-dependent boundary conditions has not been examined.

Driven by the purpose of filling the gap in the knowledge base, we dedicate this paper to the study of global dynamics of large-data classical solutions to \eqref{TOS1} subject to time-dependent Dirichlet type boundary conditions. Our major task is to identify a set of conditions on the dynamic boundary data, under which solutions starting in the $H^2$-space are globally well-posed and stabilize in the long run.  

The remainder of the paper is organized as follows: First, we state the main analytical results in Section \ref{results}. Then, we prove the main results in Sections \ref{PoT1} and  \ref{PoT2}, respectively, by developing \textit{a priori} estimates. We conclude the paper with numerical simulations that verify the analytical results and moreover, show that the assumptions on the boundary data are necessary for the results to hold.

\section{Statement of Results}\label{results}

To simplify the presentation, we take $\chi=D=1$, since the specific values of the parameters do not affect our qualitative analysis. Also, without loss of generality, we take the spatial interval as the unit interval, i.e., $(0,1)$. Moreover, recall that we concentrate upon the case of $\chi\mu>0$. Under these circumstances, the model \eqref{TOS1} is written as
\begin{subequations}\label{1.1}
\begin{alignat}{2}
u_t - (u v)_x &= u_{xx},\label{1.1a}\\
v_t - u_x &= \varepsilon v_{xx} + \varepsilon (v^2)_x.\label{1.1b}
\end{alignat}
\end{subequations}
When $\varepsilon=0$, system \eqref{1.1} becomes 
\begin{subequations}\label{1.2}
\begin{alignat}{2}
u_t-(uv)_x &=u_{xx}, \label{1.2a}\\
v_t-u_x &= 0. \label{1.2b} 
\end{alignat}
\end{subequations}
Both \eqref{1.1} and \eqref{1.2} are subject to the initial condition 
\begin{equation}\label{1.3}
(u,v)(x,0)=(u_0,v_0)(x), \quad x\in [0,1].
\end{equation}
When $\varepsilon>0$, \eqref{1.1} is supplemented with the boundary conditions:
\begin{subequations}\label{1.4}
\begin{alignat}{2}
&u(0,t)=\alpha_1(t),\quad u(1,t)=\alpha_2(t), \quad t \ge0, \label{1.4a}\\
&v(0,t)=\beta_1(t),\quad v(1,t)=\beta_2(t), \quad\ t \ge0. \label{1.4b}
\end{alignat}
\end{subequations}
When $\varepsilon=0$, only the function $u$ needs supplementary information from the boundary:
\begin{equation}\label{1.5}
\begin{aligned}
u(0,t)=\alpha_1(t),\quad u(1,t)=\alpha_2(t), \quad t \ge0,
\end{aligned}
\end{equation}
while the function $v$ does not, since otherwise the whole system would be over-determined.

With the initial and boundary conditions at our disposal, we are now ready to present the main results of this paper. The first one is concerned with the global nonlinear stability of large-data classical solutions to the diffusive model under unmatched dynamic boundary conditions for $v$. This upgrades the result of \cite{1d5}, where the case of matched boundary data is analyzed.

\begin{theorem}\label{thm1}
Consider the initial-boundary value problem \eqref{1.1}, \eqref{1.3}, \eqref{1.4}. Suppose the initial data satisfy $u_0>0$, $(u_0,v_0) \in [H^2((0,1))]^2$ and are compatible with the boundary conditions. Assume $\alpha_1=\alpha_2=\alpha$, $\beta_1$ and $\beta_2$ are smooth functions on $[0,\infty)$ satisfying  
\begin{align}
&\bullet\quad \alpha(t)\ge \underline{\alpha} >0 ,\quad \forall\,t\ge0, \quad \text{and}\quad \alpha' \in W^{1,1}(\mathbb{R}_+), \label{aa}\\
&\bullet\quad\beta_1-\beta_2 \in L^1(\mathbb{R}_+) \quad \text{and} \quad \beta_1',\beta_2' \in W^{1,1}(\mathbb{R}_+), \label{ab}
\end{align}
where $\underline{\alpha}$ is a constant. Then there exists a unique solution to the IBVP such that 
$$
\|\tilde{u}(t)\|^2_{H^2((0,1))}+\|\tilde{v}(t)\|^2_{H^2((0,1))} + \int_0^t \big(\|\tilde{u}(\tau)\|^2_{H^3((0,1))}+\|\tilde{v}(\tau)\|^2_{H^3((0,1))}\big)\mathrm{d}\tau \le C,
$$
where $\tilde{u}(x,t) = u(x,t) - \alpha(t)$, $\tilde{v}(x,t) = v(x,t) - [\beta_2(t)-\beta_1(t)]x + \beta_1(t)$, and the constant $C>0$ is independent of $t$. Moreover, the solution has the large-time behavior:
$$
\|\tilde{u}(t)\|_{H^2((0,1))}+\|\tilde{v}(t)\|_{H^2((0,1))} \to 0\quad\text{as}\quad t\to\infty.
$$ 
\end{theorem}

The second theorem addresses the dynamics of large-data classical solutions to the non-chemically-diffusive model subject to matched dynamic boundary conditions for $u$, which has not been documented in the literature. This is a generalization of the result of \cite{1d2}, where the global stability of constant Dirichlet type boundary condition is established.

\begin{theorem}\label{thm2}
Consider the initial-boundary value problem \eqref{1.2}, \eqref{1.3}, \eqref{1.5}. Suppose the initial data satisfy $u_0>0$, $(u_0,v_0) \in [H^2((0,1))]^2$ and are compatible with the boundary conditions. Assume $\alpha_1=\alpha_2=\alpha$ is a smooth function on $[0,\infty)$ satisfying
\begin{align*}
\alpha(t) \ge \underline{\alpha} >0,\quad \forall\,t\ge0, \quad \text{and}\quad \alpha' \in W^{1,1}(\mathbb{R}_+),
\end{align*}
where $\underline{\alpha}$ is a constant. Then there exists a unique solution to the IBVP such that 
$$
\|\tilde{u}(t)\|^2_{H^2((0,1))}+\|\tilde{v}(t)\|^2_{H^2((0,1))} + \int_0^t \big(\|\tilde{u}(\tau)\|^2_{H^3((0,1))}+\|\tilde{v}(\tau)\|^2_{H^2((0,1))}\big)\mathrm{d}\tau \le C,
$$
where $\tilde{u}(x,t) = u(x,t) - \alpha(t)$, $\tilde{v}(x,t) =v(x,t) - \int_0^1 v_0(x)\dx$, and the constant $C>0$ is independent of $t$. Moreover, the solution has the large-time behavior:
$$
\|\tilde{u}(t)\|_{H^2((0,1))}+\|\tilde{v}(t)\|_{H^2((0,1))} \to 0\quad\text{as}\quad t\to\infty.
$$ 
\end{theorem}

We have several remarks concerning Theorems \ref{thm1}--\ref{thm2}.

\begin{remark}\label{rem1}
The assumptions in Theorem \ref{thm1} imply the difference between $\beta_1(t)$ and $\beta_2(t)$ will converge to zero as $t\to\infty$, i.e., the unmatched boundary data will eventually match. However, based on the assumptions we see that the boundary functions are not necessarily equal to each other at any finite time. This generalizes all of the previous results for the model under matched boundary conditions. We also expect that our results will help provide useful information for the understanding of more realistic situations involving logarithmically sensitive chemotactic movements.   
\end{remark}

\begin{remark}\label{rem2}
Since the boundary functions are smooth on $[0,\infty)$ and their first order derivatives belong to $W^{1,1}(\mathbb{R}_+)$, it follows from the Fundamental Theorem of Calculus that the functions themselves alongside their first order derivatives are uniformly bounded with respect to $t$. Such information will be frequently utilized in the proof of the theorems.
\end{remark}

\begin{remark}\label{rem3}
It is an intriguing question to ask whether the matched boundary data of $u$ can be relaxed, i.e., $\alpha_1(t)\neq \alpha_2(t)$. In this case, the reference profile interpolating the boundary values becomes $\alpha_1(t) + x[\alpha_2(t)-\alpha_1(t)]$. Unfortunately, the $x$-dependence of the profile picks up additional nonlinearities when implementing the entropy estimate, see \eqref{e3}, which can not be handled by using the approach in this paper, due to the sub-quadraticity of the relative entropy. We leave the investigation for the future. 
\end{remark}

\begin{remark}
Theorem \ref{thm1} suggests that the large-time behavior of \eqref{1.1} is determined by its boundary data, while its initial information is gradually lost as time evolves. On the other hand, Theorem \ref{thm2} indicates that both its initial and boundary information are carried by \eqref{1.2} in the long run. It is interesting to investigate whether there are mechanisms that can drive the solutions to other stationary solutions rather than those determined by the initial and/or boundary data. Among various types of modeling structures, the logistic growth may give us a definite answer. However, the analysis in this paper can not  be directly carried over to the model with logistic growth, due to the quadratic nonlinearity. We will report such a result in a forthcoming paper. 
\end{remark}

\begin{remark}
Systems \eqref{1.1} and \eqref{1.2} themselves have deep mathematical interests as they serve as prototypes of general parabolic/hyperbolic balance laws. The nature of possible non-uniform dissipativity (i.e., $\varepsilon=0$) coupled with nonlinear flux functions, together with the non-triviality of the dynamic boundary data in this type of problems presents significant challenges in mathematical analysis. The analysis of \eqref{1.1} and \eqref{1.2} helps to shed light on how to advance fundamental research of parabolic/hyperbolic balance laws in related topics, and we expect the study in this paper to be helpful to other parabolic/hyperbolic balance laws. 
\end{remark}

We prove Theorem \ref{thm1} and Theorem \ref{thm2} in Section \ref{PoT1} and Section \ref{PoT2}, respectively. The integrability conditions of the boundary data are thoroughly examined to fit into the framework previously established for the case of matched boundary conditions. The proof of Theorem \ref{thm1} takes advantage of the fully dissipative structure of \eqref{1.1}, which results in energy estimates depending on the reciprocal of the chemical diffusion coefficient. For this reason, the arguments can not be carried over to \eqref{1.2}. Instead, Theorem \ref{thm2} is proven by deriving a nonlinear damping equation for the spatial derivative of $v$. It should be mentioned that the approach for analyzing \eqref{1.2} can not be utilized for \eqref{1.1}, due to the lack of information of higher order spatial derivatives of the solution to the latter, causing integration-by-parts to be unaccessible. Thus the energy methods for studying \eqref{1.1} and \eqref{1.2} are mutually exclusive.

Lastly, for notational convenience, throughout the rest of the paper, we use $\|\cdot\|$, $\|\cdot\|_{H^s}$ and $\|\cdot\|_\infty$ to denote the standard norms $\|\cdot\|_{L^2((0,1))}$, $\|\cdot\|_{H^s((0,1))}$ and $\|\cdot\|_{L^\infty((0,1))}$, respectively. Moreover, we use $C$ to denote a generic constant which is independent of time, but may depend on the parameter $\varepsilon$ and initial and/or boundary data. The value of the constant may vary line by line according to the context.

\section{Proof of Theorem \ref{thm1}}\label{PoT1}

The proof of Theorem \ref{thm1} is divided into five steps contained in a series of subsections. First of all, the local well-posedness of the IBVP under the assumptions of Theorem \ref{thm1} can be established via standard approaches, such as mollification, Galerkin approximation, energy estimate, contraction mapping principle and compactness argument. We omit most of the standard technical details for brevity, while focus on deriving the {\it a priori} estimates of the local solution, in order to extend it to a global one. Moreover, the {\it a priori} estimates can be justified by standard means, which are indeed carried out on the local smooth approximate solutions obtained from the mollified initial data and contraction mapping. We begin with an estimate based on the relative entropy-entropy flux pair associated with the IBVP.

\subsection{Entropy Estimate}

\begin{lemma}\label{lem1}
Under the assumptions of Theorem $\ref{thm1}$, there exists a constant $C>0$ which is independent of $t$, such that
\begin{equation*}
\begin{aligned}
&E(u,\alpha)(t)+\|(v-\beta)(t)\|^2+\int_0^t\int_0^1\frac{(u_x)^2}{u}\mathrm{d}x\mathrm{d}\tau+\int_0^t \varepsilon\|(v_x-\beta_x)(\tau)\|^2\mathrm{d}\tau \le C,
\end{aligned}
\end{equation*}
where $\beta(x,t)=[\beta_2(t)-\beta_1(t)]x + \beta_1(t)$ and 
\begin{equation*}
E(u,\alpha)\equiv \int_0^1 \left[(u\ln u-u)-(\alpha\ln\alpha-\alpha)-(u-\alpha)\ln\alpha\right]\mathrm{d}x\ge0
\end{equation*}
denotes the relative entropy.
\end{lemma}

\begin{proof} First of all, we note that according to the local well-posedness theory and the assumptions on the initial and boundary data, especially $u_0>0$ and $\alpha\ge \underline{\alpha}>0$, the function $u$ is positive within the life span of the local solution, i.e, $[0,T^*)$ for some $T^*>0$. The subsequent estimates are derived within such a time window. It will be shown that the estimates are indeed independent of $t$. Then the global well-posedness follows from the uniform estimates and standard continuation argument.

{\bf Step 1.} By a direct calculation, we can show that
\begin{align}\label{e2}
(u\ln u-u)_t-(\alpha\ln\alpha-\alpha)_t-[(u-\alpha)\ln\alpha]_t 
=(\ln u-\ln\alpha)u_t-(u-\alpha)\frac{\alpha'}{\alpha}.
\end{align}
Using equation \eqref{1.1a} and noting {$\alpha$} depends only on {$t$}, we deduce:
\begin{align}\label{e3}
(\ln u-\ln\alpha)u_t
=[(\ln u-\ln\alpha)uv]_x+[(\ln u-\ln\alpha)u_x]_x-v\,u_x-\frac{(u_x)^2}{u}.
\end{align}
Substituting \eqref{e3} into \eqref{e2} gives us
\begin{align}\label{e4}
&(u\ln u-u)_t-(\alpha\ln\alpha-\alpha)_t-[(u-\alpha)\ln\alpha]_t \notag\\
=&[(\ln u-\ln\alpha)uv]_x+[(\ln u-\ln\alpha)u_x]_x-v u_x-\frac{(u_x)^2}{u}-(u-\alpha)\frac{\alpha'}{\alpha}.
\end{align}
Integrating \eqref{e4} over {$[0,1]$} and using the boundary conditions, we have
\begin{equation}\label{e5}
\begin{aligned}
\frac{\mathrm{d}}{\mathrm{d}t}E(u,\alpha)+\int_0^1\frac{(u_x)^2}{u}\mathrm{d}x = -\int_0^1 v u_x\mathrm{d}x-\int_0^1(u-\alpha)\frac{\alpha'}{\alpha}\mathrm{d}x.
\end{aligned}
\end{equation}
Let $\beta(x,t)=[\beta_2(t)-\beta_1(t)]x + \beta_1(t)$. Then we derive from equation \eqref{1.1b} that
\begin{align}\label{e6}
&(v-\beta)_t-u_x\notag\\
=&\varepsilon(v-\beta)_{xx}-2\varepsilon (v-\beta)(v-\beta)_x -2\varepsilon (v-\beta)\beta_x-2\varepsilon \beta (v-\beta)_x - 2\varepsilon \beta\beta_x-\beta_t.
\end{align}
Taking {$L^2$} inner product of \eqref{e6} with {$v-\beta$} and using the boundary conditions, we obtain
\begin{align}\label{e7}
&\frac12\frac{\mathrm{d}}{\mathrm{d}t} \|v-\beta\|^2+ \varepsilon\|(v-\beta)_x\|^2\notag\\
= &\int_0^1(v-\beta)u_x\mathrm{d}x - \varepsilon\int_0^1(v-\beta)^2\beta_x\mathrm{d}x -2\varepsilon\int_0^1\beta\beta_x(v-\beta)\mathrm{d}x -\int_0^1(v-\beta)\beta_t\mathrm{d}x.
\end{align}
Note since $\alpha$ is independent of $x$, it holds that
\begin{equation*}
\begin{aligned}
\int_0^1(v-\beta)u_x\mathrm{d}x 
=\int_0^1 v u_x\mathrm{d}x + \int_0^1 \beta_x(u-\alpha)\mathrm{d}x.
\end{aligned}
\end{equation*}
We then get from \eqref{e7} that
\begin{align}\label{e8}
&\frac12\frac{\mathrm{d}}{\mathrm{d}t} \|v-\beta\|^2+ \varepsilon\|(v-\beta)_x\|^2 \notag\\
=&\int_0^1 v u_x\mathrm{d}x+ \int_0^1 \beta_x(u-\alpha)\mathrm{d}x - \varepsilon\int_0^1(v-\beta)^2\beta_x\mathrm{d}x -2\varepsilon\int_0^1\beta\beta_x(v-\beta)\mathrm{d}x\notag\\
& -\int_0^1(v-\beta)\beta_t\mathrm{d}x.
\end{align}
Adding \eqref{e8} and \eqref{e5}, we can show that
\begin{align}\label{e9}
&\frac{\mathrm{d}}{\mathrm{d}t}\Big(E(u,\alpha)+\frac12\|v-\beta\|^2\Big)+\int_0^1\frac{(u_x)^2}{u}\mathrm{d}x+ \|(v-\beta)_x\|^2 \notag\\
=&\int_0^1(u-\alpha)\big(\beta_x-\frac{\alpha'}{\alpha}\big)\mathrm{d}x - \varepsilon \int_0^1(v-\beta)^2\beta_x\mathrm{d}x - \int_0^1 (2 \varepsilon\beta\beta_x+\beta_t)(v-\beta)\mathrm{d}x \notag \\
\le&\Big(\frac{|\alpha'|}{\alpha}+|\beta_x|\Big)\int_0^1|u-\alpha|\mathrm{d}x+\varepsilon|\beta_x|\,\|\beta-v\|^2+\int_0^1\big(2\varepsilon|\beta\beta_x|+|\beta_t|\big)|v-\beta|\mathrm{d}x.
\end{align}

{\bf Step 2.} Define
$$
F_\alpha(u)\equiv (u\ln u-u)-(\alpha\ln\alpha-\alpha)-(u-\alpha)\ln\alpha+(e-1)\alpha-u.
$$
Then it can be readily checked that $F_\alpha(e\alpha)=0$, $F_\alpha'(e\alpha)=0$, and $F_\alpha''(u)=\frac{1}{u}\ge0$ for $u\ge0$. These imply $F_\alpha(u)\ge0$ for $u\ge0$. Hence,
$$
0\le u \le (u\ln u-u)-(\alpha\ln\alpha-\alpha)-(u-\alpha)\ln\alpha+(e-1)\alpha,
$$
and  
\begin{equation}\label{L1}
0\le \int_0^1 u(x,t)\mathrm{d}x \le E(u,\alpha)+(e-1)\alpha.
\end{equation}
By triangle inequality, we have
\begin{align}\label{L1a}
\int_0^1|u-\alpha|\mathrm{d}x \le E(u,\alpha)+e\alpha.
\end{align}
Applying \eqref{L1a} to the right-hand side of \eqref{e9}, we can show that
\begin{align*}
\bullet&\quad \Big(\frac{|\alpha'|}{\alpha}+|\beta_x|\Big)\int_0^1|u-\alpha|\mathrm{d}x \le C\big(|\alpha'| + |\beta_1 - \beta_2|\big)\big[E(u,\alpha)+1\big],\\
\bullet&\quad \varepsilon |\beta_x|\,\|\beta-v\|^2 \le \varepsilon |\beta_1-\beta_2|\,\|\beta-v\|^2,\\
\bullet&\quad\int_0^1\big(2\varepsilon |\beta\beta_x|+|\beta_t|\big)|v-\beta|\mathrm{d}x \le C\big(\varepsilon |\beta_1-\beta_2|+|\beta_1'|+|\beta_2'|\big)\big(\|v-\beta\|^2+1\big),
\end{align*}
where we used the uniform boundedness of $\beta$ and the Cauchy-Schwarz inequality. Using the above estimates, we update \eqref{e9} as
\begin{align}\label{e10}
&\frac{\mathrm{d}}{\mathrm{d}t}\Big(E(u,\alpha)+\frac12\|v-\beta\|^2+1\Big)+\int_0^1\frac{(u_x)^2}{u}\mathrm{d}x+ \varepsilon \|(v-\beta)_x\|^2 \notag\\
\le\ &C\big(|\alpha'| + \varepsilon |\beta_1-\beta_2| + |\beta_1'|+ |\beta_2'|\big)\Big(E(u,\alpha)+\frac12\|v-\beta\|^2+1\Big),
\end{align}
where the constant is independent of $t$. Applying Gr\"onwall's inequality to \eqref{e10} gives us
\begin{align}\label{e11}
&E(u,\alpha)(t)+\frac12\|(v-\beta)(t)\|^2+1\notag \\
\le\ &\exp\Big\{C\int_0^t\big(|\alpha'| + \varepsilon |\beta_1-\beta_2| + |\beta_1'|+ |\beta_2'|\big)\mathrm{d}\tau\Big\}\Big(E(u_0,\alpha_0)+\frac12\|v_0-\beta_0\|^2+1\Big).
\end{align}
Using the assumptions in Theorem \ref{thm1}, we deduce from \eqref{e11} that
\begin{equation}\label{e13}
E(u,\alpha)(t)+\frac12\|(v-\beta)(t)\|^2\le C,
\end{equation}
where the constant is independent of $t$. Substituting \eqref{e13} into \eqref{e10}, then integrating the resulting inequality with respect to $t$, we have in particular,
\begin{equation}\label{e14}
\int_0^t\int_0^1\frac{(u_x)^2}{u}\mathrm{d}x\mathrm{d}\tau+\int_0^t \varepsilon \|(v_x-\beta_x)(\tau)\|^2\mathrm{d}\tau \le C,
\end{equation}
where the constant is independent of $t$. This, along with \eqref{e13}, completes the entropy estimate and hence the proof of Lemma \ref{lem1}.
\end{proof}

\subsection{$L^\infty_tL^2_x$--$L^2_tH^1_x$--Estimates} We now switch to standard $L^2$-based energy estimates. To facilitate our asymptotic analysis, we define
$$
\tilde{u}\equiv u-\alpha\quad \text{and} \quad \tilde{v}\equiv v-\beta,
$$
where $(u, v)$ satisfies \eqref{1.1} and $\beta(x,t)=[\beta_2(t)-\beta_1(t)]x+\beta_1(t)$. Then {$(\tilde{u},\tilde{v})$} satisfies
\begin{subequations}\label{e16}
\begin{alignat}{4}
\tilde{u}_t-(\tilde{u}\tilde{v})_x&=\tilde{u}_{xx} +\alpha \tilde{v}_x + \beta\tilde{u}_x + \beta_x\tilde{u}+\alpha\beta_x-\alpha', \label{e16a}\\
\tilde{v}_t-\tilde{u}_x&= \varepsilon \tilde{v}_{xx}-2 \varepsilon  \tilde{v}\tilde{v}_x -2\varepsilon \beta\tilde{v}_x -2\varepsilon \tilde{v}\beta_x -2\varepsilon \beta\beta_x-\beta_t, \label{e16b}\\
(\tilde{u},\tilde{v})(x,0)&=(u_0(x)-\alpha(0),v_0(x)-\beta(x,0)), \label{e16c}\\
\tilde{u}(0,t)&=\tilde{u}(1,t)=0,\ \ \ \tilde{v}(0,t)=\tilde{v}(1,t)=0. \label{e16d}
\end{alignat}
\end{subequations}
Using Lemma \ref{lem1}, we can show the following:

\begin{lemma}\label{lem2}
Under the assumptions of Theorem $\ref{thm1}$, there exists a constant $C>0$ which is independent of $t$, such that
\begin{equation*}\label{e120}
\begin{aligned}
\|\tilde{u}(t)\|^2+\|\tilde{v}(t)\|^2+\int_0^t\|\tilde{u}_x(\tau)\|^2\mathrm{d}\tau \le C.
\end{aligned}
\end{equation*}
\end{lemma}

\begin{proof}
{\bf Step 1.} Taking $L^2$ inner product of \eqref{e16a} with $\tilde{u}$, we have
\begin{equation}\label{e17}
\frac12\frac{\mathrm{d}}{\mathrm{d}t}\|\tilde{u}\|^2+\|\tilde{u}_x\|^2=-\int_0^1 \tilde{u}\tilde{v}\tilde{u}_x\mathrm{d}x+\alpha\int_0^1 \tilde{u}\tilde{v}_x\mathrm{d}x +(\alpha\beta_x-\alpha')\int_0^1\tilde{u}\mathrm{d}x +\frac{\beta_x}{2}\|\tilde{u}\|^2.
\end{equation}
Taking {$L^2$} inner product of \eqref{e16b} with $\tilde{v}$ yields
\begin{equation}\label{e18}
\frac12\frac{\mathrm{d}}{\mathrm{d}t}\|\tilde{v}\|^2+\varepsilon \|\tilde{v}_x\|^2=\int_0^1 \tilde{v}\tilde{u}_x\mathrm{d}x - \varepsilon \beta_x\|\tilde{v}\|^2 -\int_0^1 (2\varepsilon \beta\beta_x+\beta_t)\tilde{v}\mathrm{d}x.
\end{equation}
Multiplying \eqref{e18} by $\alpha$, we obtain
\begin{equation}\label{e19}
\begin{aligned}
\frac12\frac{\mathrm{d}}{\mathrm{d}t}\left(\alpha\|\tilde{v}\|^2\right)+\alpha\varepsilon \|\tilde{v}_x\|^2
&=-\alpha\int_0^1 \tilde{u}\tilde{v}_x\mathrm{d}x +\big(\frac{\alpha'}{2}-\alpha\varepsilon \beta_x\big)\|\tilde{v}\|^2 - \alpha\int_0^1 (2\varepsilon \beta\beta_x+\beta_t) \tilde{v}\mathrm{d}x.
\end{aligned}
\end{equation}
Adding \eqref{e19} to \eqref{e17} gives us
\begin{align}\label{e20}
&\frac12\frac{\mathrm{d}}{\mathrm{d}t}\left(\|\tilde{u}\|^2+\alpha \|\tilde{v}\|^2\right)+\|\tilde{u}_x\|^2+\alpha \varepsilon\|\tilde{v}_x\|^2 \notag\\
=&\underbrace{-\int_0^1 \tilde{u}\tilde{v}\tilde{u}_x\mathrm{d}x}_{\equiv I_1}+\underbrace{(\alpha\beta_x-\alpha')\int_0^1\tilde{u}\mathrm{d}x}_{\equiv I_2} +\underbrace{\frac{\beta_x}{2}\|\tilde{u}\|^2}_{\equiv I_3} +\underbrace{\big(\frac{\alpha'}{2}-\alpha\varepsilon \beta_x\big)\|\tilde{v}\|^2}_{\equiv I_4}\notag\\
& \underbrace{-\alpha\int_0^1 (2\varepsilon \beta\beta_x+\beta_t)\tilde{v}\mathrm{d}x}_{\equiv I_5}.
\end{align}

{\bf Step 2.} To estimate $I_1$, we note that 
\begin{align}\label{e21}
|I_1| \le \frac12\|\tilde{u}\|_{L^\infty}^2\|\tilde{v}\|^2 + \frac12\|\tilde{u}_x\|^2 \le C\|\tilde{u}\|_{L^\infty}^2 + \frac12\|\tilde{u}_x\|^2,
\end{align}
where we used \eqref{e13}. Since $\tilde{u}(0,t)=0$, for any $x\in [0,1]$, it holds that
\begin{align}\label{e22}
|\tilde{u}(x,t)|^2=\Big|\int_0^x \tilde{u}_y(y,t)\mathrm{d}y\Big|^2 &\le \Big(\int_0^1|\tilde{u}_x|\mathrm{d}x \Big)^2 \notag\\
&\le \Big(\int_0^1u(x,t)\mathrm{d}x\Big)\Big(\int_0^1\frac{(\tilde{u}_x)^2}{u}\mathrm{d}x\Big),
\end{align}
where $u$ denotes the solution to \eqref{1.1}. Since $E(u,\alpha)$ and $\alpha$ are uniformly bounded with respect to $t$ (see \eqref{e13}), we obtain from \eqref{L1} and \eqref{e22} that
\begin{align}\label{e23}
\|\tilde{u}(t)\|_{L^\infty}^2\le C\int_0^1\frac{(\tilde{u}_x)^2}{u}\mathrm{d}x.
\end{align}
Substituting \eqref{e23} into \eqref{e21} gives us
\begin{equation}\label{e23a}
\left|I_1\right|\le C\int_0^1\frac{(u_x)^2}{u}\mathrm{d}x+\frac12\|\tilde{u}_x\|^2,
\end{equation}
where $C$ is independent of $t$. The remaining terms are estimated as
$$
|I_2|\le C\big(|\beta_1-\beta_2|+|\alpha'| \big)\big(\|\tilde{u}\|^2+1\big),
$$
and
$$
\begin{aligned}
|I_5| \le  C\big(\varepsilon|\beta_1-\beta_2|+|\beta_1'|+|\beta_2'|\big)\big(\|\tilde{v}\|^2+1\big).
\end{aligned}
$$
Then we update \eqref{e20} as
\begin{align}\label{e23a}
&\ \frac12\frac{\mathrm{d}}{\mathrm{d}t}\left(\|\tilde{u}\|^2+\alpha\|\tilde{v}\|^2+1\right)+\frac12\|\tilde{u}_x\|^2+\alpha\varepsilon\|\tilde{v}_x\|^2\notag\\
\le&\ C\big[|\alpha'| + (1+\varepsilon) |\beta_1-\beta_2| + |\beta_1'|+|\beta_2'| \big]\left(\|\tilde{u}\|^2+\alpha\|\tilde{v}\|^2+1\right) + C\int_0^1\frac{(u_x)^2}{u}\mathrm{d}x,
\end{align}
where the constants are independent of $t$. Applying Gr\"onwall's inequality to \eqref{e23a} and using \eqref{e14} and the assumptions in Theorem \ref{thm1}, we find that
\begin{equation}\label{e24}
\begin{aligned}
\|\tilde{u}(t)\|^2+\alpha(t)\|\tilde{v}(t)\|^2\le C,
\end{aligned}
\end{equation}
where the constant is independent of $t$. Substituting \eqref{e24} into \eqref{e23a}, then integrating the resulting inequality with respect to $t$, we can show that
\begin{equation}\label{e25}
\begin{aligned}
\int_0^t\|\tilde{u}_x(\tau)\|^2\mathrm{d}\tau\le C,
\end{aligned}
\end{equation}
where the constant is independent of $t$. We conclude the proof by noticing $0< \underline{\alpha}\le \alpha(t)$. This completes the proof of Lemma \ref{lem2}.
\end{proof}

\subsection{$L^\infty_tH^1_x$--$L^2_tH^2_x$--Estimates}
 \begin{lemma}\label{lem3}
 Under the assumptions of Theorem $\ref{thm1}$, there exists a constant $C>0$ which is independent of $t$, such that
\begin{equation*}
\begin{aligned}
\|\tilde{u}_x(t)\|^2+\|\tilde{v}_x(t)\|^2+\int_0^t\left(\|\tilde{u}_{xx}(\tau)\|^2+\varepsilon\|\tilde{v}_{xx}(\tau)\|^2\right)\mathrm{d}\tau
\le C.
\end{aligned}
\end{equation*}
\end{lemma}

\begin{proof}
Taking $L^2$ inner products of \eqref{e16a} with $-\tilde{u}_{xx}$ and \eqref{e16b} with $-\tilde{v}_{xx}$, respectively, then adding the results, we obtain
\begin{align}\label{e26}
&\frac12\frac{\mathrm{d}}{\mathrm{d}t}\left(\|\tilde{u}_x\|^2+\|\tilde{v}_x\|^2\right) +\|\tilde{u}_{xx}\|^2 + \varepsilon\|\tilde{v}_{xx}\|^2 \notag\\
=&\underbrace{-\int_0^1\left(\tilde{v}\tilde{u}_x+\tilde{u}\tilde{v}_x+\alpha \tilde{v}_x+\beta\tilde{u}_x+\beta_x\tilde{u}\right)\,\tilde{u}_{xx}\mathrm{d}x}_{\equiv J_1} \underbrace{-(\alpha\beta_x-\alpha')\int_0^1\tilde{u}_{xx}\mathrm{d}x}_{\equiv J_2} +\underbrace{2\varepsilon\int_0^1 \tilde{v}\tilde{v}_x\tilde{v}_{xx}\mathrm{d}x}_{\equiv J_3} \notag\\
&+\underbrace{2\varepsilon\int_0^1 \beta\tilde{v}_x\tilde{v}_{xx}\mathrm{d}x}_{\equiv J_4} +\underbrace{2\varepsilon\beta_x\int_0^1 \tilde{v}\tilde{v}_{xx}\mathrm{d}x}_{\equiv J_5} \underbrace{- \int_0^1 \tilde{u}_x\tilde{v}_{xx}\mathrm{d}x}_{\equiv J_6}+\underbrace{\int_0^1 (2\varepsilon\beta\beta_x+\beta_t)\tilde{v}_{xx}\mathrm{d}x}_{\equiv J_7}.
\end{align}
Since $\beta=(\beta_2-\beta_1)x+\beta_1$ and $\alpha,\beta_1,\beta_2$ are uniformly bounded (see Remark \ref{rem2}), using the Cauchy-Schwarz, Sobolev and Poincar\'e inequalities, we can show that 
$$
\begin{aligned}
|J_1| &\le \frac14\|\tilde{u}_{xx}\|^2+C\left(\|\tilde{v}\|_{L^\infty}^2\|\tilde{u}_x\|^2+\|\tilde{u}\|^2_{L^\infty}\|\tilde{v}_x\|^2+\|\tilde{v}_x\|^2+\|\tilde{u}_x\|^2+\|\tilde{u}\|^2\right)\\
&\le \frac14\|\tilde{u}_{xx}\|^2 + C\left(\|\tilde{u}_x\|^2\|\tilde{v}_x\|^2+\|\tilde{v}_x\|^2+\|\tilde{u}_x\|^2\right).
\end{aligned}
$$
Since $\alpha'(t)$ is uniformly bounded (see Remark \ref{rem2}), we estimate $J_2$ as:
\begin{align*}
|J_2| &\le \frac14\|\tilde{u}_{xx}\|^2+ C(|\beta_1-\beta_2|^2+|\alpha'|^2) \le \frac14\|\tilde{u}_{xx}\|^2+ C(|\beta_1-\beta_2|+|\alpha'|).
\end{align*}
Similar to the estimate of $J_1$, we can show that
$$
|J_3| \le \frac{\varepsilon}{4}\|\tilde{v}_{xx}\|^2+C\varepsilon\|\tilde{v}_x\|^2\|\tilde{v}_x\|^2.
$$
By the uniform boundedness of the boundary data and Poincar\'e's inequality, we estimate the remaining terms as:
$$
|J_4|+|J_5|+|J_6|+|J_7| \le \frac{\varepsilon}{4}\|\tilde{v}_{xx}\|^2+C(\varepsilon,\varepsilon^{-1})\big(\|\tilde{v}_x\|^2+\|\tilde{u}_x\|^2+|\beta_1-\beta_2|+|\beta_1'|+|\beta_2'|\big).
$$
We remark that the constant $C(\varepsilon,\varepsilon^{-1})$ results from applying the Cauchy-Schwarz inequality. Substituting the above estimates into \eqref{e26} gives us
\begin{align}\label{e27}
&\ \frac12\frac{\mathrm{d}}{\mathrm{d}t}\left(\|\tilde{u}_x\|^2+\|\tilde{v}_x\|^2\right)+\frac12\|\tilde{u}_{xx}\|^2+\frac{\varepsilon}{2}\|\tilde{v}_{xx}\|^2 \notag\\
\le &\ C\|\tilde{v}_x\|^2(\|\tilde{u}_x\|^2+\varepsilon\|\tilde{v}_x\|^2)+C(\varepsilon,\varepsilon^{-1})\big(\|\tilde{v}_x\|^2+\|\tilde{u}_x\|^2+|\alpha'|+|\beta_1-\beta_2|+|\beta_1'|+|\beta_2'|\big).
\end{align}
Applying Gr\"onwall's inequality to \eqref{e27} and using \eqref{e14}, \eqref{e25} and the assumptions in Theorem \ref{thm1}, we can show that
\begin{equation}\label{e28}
\begin{aligned}
\|\tilde{u}_x(t)\|^2+\|\tilde{v}_x(t)\|^2+\int_0^t\left(\|\tilde{u}_{xx}(\tau)\|^2+\varepsilon\|\tilde{v}_{xx}(\tau)\|^2\right)\mathrm{d}\tau \le C(\varepsilon,\varepsilon^{-1}),
\end{aligned}
\end{equation}
where the constant is independent of $t$. This completes the proof of Lemma \ref{lem3}.
\end{proof}

\subsection{$L^\infty_tH^2_x$--$L^2_tH^3_x$--Estimates}
 \begin{lemma}\label{lem4}
 Under the assumptions of Theorem $\ref{thm1}$, there exists a constant $C>0$ which is independent of $t$, such that
\begin{equation*}
\begin{aligned}
\|\tilde{u}_{xx}(t)\|^2+\|\tilde{v}_{xx}(t)\|^2+\int_0^t\left(\|\tilde{u}_{xxx}\|^2+ \varepsilon \|\tilde{v}_{xxx}\|^2\right)\mathrm{d}\tau
\le C.
\end{aligned}
\end{equation*}
\end{lemma}

\begin{proof} Since the information of the higher order spatial derivatives of the solution is unknown at the boundary points, the usual procedure (differentiating with respect to $x$) for estimating the $L^\infty_tH^2_x$ and $L^2_tH^3_x$ norms of the solution can not be directly implemented here. To circumvent such a technical obstruction, we turn to the estimation of the temporal derivatives of the solution, then utilize the equations to recover the spatial derivatives.

{\bf Step 1.} Taking $\partial_t$ of \eqref{e16a} and \eqref{e16b}, we obtain
\begin{subequations}\label{e29}
\begin{alignat}{3}
\tilde{u}_{tt}-(\tilde{u}\tilde{v})_{xt} &=\tilde{u}_{xxt} +\alpha' \tilde{v}_x+\alpha \tilde{v}_{xt}+\beta_t\tilde{u}_x+\beta\tilde{u}_{xt}+\beta_{xt}\tilde{u}+\beta_x\tilde{u}_t+\alpha'\beta_x-\alpha\beta_{xt}-\alpha'', \label{e29a}\\
\tilde{v}_{tt}-\tilde{u}_{xt} &=\varepsilon \tilde{v}_{xxt} - 2 \varepsilon \tilde{v}_t\tilde{v}_x - 2 \varepsilon \tilde{v}\tilde{v}_{xt} - 2 \varepsilon \beta_t\tilde{v}_x - 2 \varepsilon \beta\tilde{v}_{xt} - 2 \varepsilon \tilde{v}_t\beta_x - 2 \varepsilon \tilde{v}\beta_{xt} - 2 \varepsilon \beta_t\beta_x \notag\\
&\qquad - 2 \varepsilon \beta\beta_{xt}-\beta_{tt}.\label{e29b}
\end{alignat}
\end{subequations}
Taking $L^2$ inner product of \eqref{e29a} with $\tilde{u}_t$, we have 
\begin{align}\label{H21}
&\frac{1}{2}\frac{\mathrm{d}}{\mathrm{d}t}\|\tilde{u}_t\|^2+\|\tilde{u}_{xt}\|^2\notag\\
=&\underbrace{-\int_0^1(\tilde{u}\tilde{v})_t\tilde{u}_{xt}\mathrm{d}x}_{\equiv K_{1}}+\underbrace{\alpha'\int_0^1\tilde{v}_x\tilde{u}_t\mathrm{d}x}_{\equiv K_{2}}\underbrace{-\alpha\int_0^1\tilde{v}_{t}\tilde{u}_{xt}\mathrm{d}x}_{\equiv K_{3}} +\underbrace{\int_0^1 \beta_t\tilde{u}_{x}\tilde{u}_t\mathrm{d}x}_{\equiv K_{4}} + \underbrace{\int_0^1 \beta\tilde{u}_{xt}\tilde{u}_t\mathrm{d}x}_{K_{5}} \notag\\
&+\underbrace{\beta_{xt}\int_0^1\tilde{u}\tilde{u}_t\mathrm{d}x}_{\equiv K_{6}}+\underbrace{\beta_x\|\tilde{u}_{t}\|^2}_{\equiv K_{7}}+\underbrace{(\alpha'\beta_x-\alpha\beta_{xt}-\alpha'')\int_0^1\tilde{u}_t\mathrm{d}x}_{\equiv K_{8}}.
\end{align}
For $K_{1}$, we can show that
$$
\begin{aligned}
|K_{1}| \le & \frac16 \|\tilde{u}_{xt}\|^2 + 3\big(\|\tilde{u}\|_{L^\infty}^2\|\tilde{v}_t\|^2 + \|\tilde{v}\|_{L^\infty}^2\|\tilde{u}_t\|^2\big)\\
\le &  \frac16 \|\tilde{u}_{xt}\|^2 + C\big(\|\tilde{u}\|_{H^1}^2\|\tilde{v}_t\|^2 + \|\tilde{v}\|_{H^1}^2\|\tilde{u}_t\|^2\big)\notag\\
\le & \frac16 \|\tilde{u}_{xt}\|^2 + C\big(\|\tilde{v}_t\|^2 + \|\tilde{u}_t\|^2\big),
\end{aligned}
$$
where we applied Lemmas \ref{lem2}--\ref{lem3}. Using the boundedness of the boundary data, we can show that
$$
|K_{2}| + |K_{3}| + |K_{4}| + |K_{5}| +|K_6| + |K_7| \le \frac13 \|\tilde{u}_{xt}\|^2 + C\big(\|\tilde{u}_x\|^2 + \|\tilde{v}_x\|^2+\|\tilde{u}_t\|^2 + \|\tilde{v}_t\|^2\big),
$$
where we also invoked Poincar\'e's inequality. Similarly, $K_{8}$ is estimated as 
$$
\begin{aligned}
|K_{8}| \le & |\alpha''|\,\|\tilde{u}_t\|^2 + C\big(|\alpha''|+\|\tilde{u}_t\|^2+|\alpha'|^2|\beta_1-\beta_2|^2+|\alpha|^2(|\beta_1'|^2+|\beta_2'|^2)\big)\\
\le & |\alpha''|\,\|\tilde{u}_t\|^2 + C\big(|\alpha''|+\|\tilde{u}_t\|^2+|\beta_1-\beta_2|+|\beta_1'|+|\beta_2'|\big).
\end{aligned}
$$
Substituting the above estimates into \eqref{H21} gives us
\begin{align}\label{H22}
\frac{1}{2}\frac{\mathrm{d}}{\mathrm{d}t}\|\tilde{u}_t\|^2+\frac12\|\tilde{u}_{xt}\|^2 \le  |\alpha''|\,\|\tilde{u}_t\|^2+C\big(&\|\tilde{u}_t\|^2+\|\tilde{u}_x\|^2+\|\tilde{v}_t\|^2+\|\tilde{v}_x\|^2 \notag\\
& +|\alpha''|+|\beta_1-\beta_2|+|\beta_1'|+|\beta_2'|\big).
\end{align}

{\bf Step 2.} Taking $L^2$ inner product of \eqref{e29b} with $\tilde{v}_t$, we obtain
\begin{align}\label{H22a}
&\frac{1}{2}\frac{\mathrm{d}}{\mathrm{d}t}\|\tilde{v}_t\|^2+ \varepsilon \|\tilde{v}_{xt}\|^2\notag\\
=&-\int_0^1\tilde{u}_t\tilde{v}_{xt}\mathrm{d}x+2 \varepsilon \int_0^1\tilde{v}\tilde{v}_t\tilde{v}_{xt}\mathrm{d}x - 2 \varepsilon \int_0^1\beta_t\tilde{v}_{x}\tilde{v}_t\mathrm{d}x - \varepsilon \beta_x\|\tilde{v}_t\|^2\notag\\
& -2 \varepsilon \beta_{xt} \int_0^1 \tilde{v}\tilde{v}_t\mathrm{d}x - 2 \varepsilon \beta_{x}\int_0^1\beta_t \tilde{v}_t\mathrm{d}x - 2 \varepsilon \beta_{xt}\int_0^1\beta \tilde{v}_t\mathrm{d}x-\int_0^1\beta_{tt}\tilde{v}_t\mathrm{d}x.
\end{align}
Using the arguments in {\bf Step 1}, we can show that
\begin{align}\label{H23}
\frac{1}{2}\frac{\mathrm{d}}{\mathrm{d}t}\|\tilde{v}_t\|^2+\frac{\varepsilon}{2}\|\tilde{v}_{xt}\|^2 \le \big(|\beta_1''|+|\beta_2''|\big)\|\tilde{v}_t\|^2+C\big(&\|\tilde{u}_t\|^2+\|\tilde{u}_x\|^2+\|\tilde{v}_t\|^2+\|\tilde{v}_x\|^2 \notag\\
& +|\beta_1-\beta_2|+|\beta_1'|+|\beta_2'|+|\beta_1''|+|\beta_2''|\big).
\end{align}
Note that according to the assumptions of Theorem \ref{thm1} and Lemmas \ref{lem1}--\ref{lem2}, the quantities on the right-hand sides of \eqref{H22} and \eqref{H23}, except $\|\tilde{u}_t\|^2$ and $\|\tilde{v}_t\|^2$, are uniformly integrable with respect to $t$. For $\|\tilde{u}_t\|^2$ and $\|\tilde{v}_t\|^2$, based on the equations \eqref{e16a} and \eqref{e16b}, we can show that
$$
\|\tilde{u}_t\|^2 \le C\big(\|\tilde{u}_x\|^2 + \|\tilde{v}_x\|^2 + \|\tilde{u}_{xx}\|^2 + |\beta_1-\beta_2| + |\alpha'|\big),
$$
and 
$$
\|\tilde{v}_t\|^2 \le C\big(\|\tilde{u}_x\|^2 + \|\tilde{v}_x\|^2 + \|\tilde{v}_{xx}\|^2 + |\beta_1-\beta_2| + |\beta_1'| + |\beta_2'|\big).
$$
Then we update \eqref{H22} and \eqref{H23} as
\begin{align}\label{H24}
\frac{1}{2}\frac{\mathrm{d}}{\mathrm{d}t}\|\tilde{u}_t\|^2+\frac12\|\tilde{u}_{xt}\|^2 \le |\alpha''|\,\|\tilde{u}_t\|^2+C\big(&\|\tilde{u}_{xx}\|^2+\|\tilde{u}_x\|^2+\|\tilde{v}_{xx}\|^2+\|\tilde{v}_x\|^2 \notag\\
& +|\alpha'|+|\alpha''|+|\beta_1-\beta_2|+|\beta_1'|+|\beta_2'|\big),
\end{align}
and 
\begin{align}\label{H25}
\frac{1}{2}\frac{\mathrm{d}}{\mathrm{d}t}\|\tilde{v}_t\|^2+\frac12\|\tilde{v}_{xt}\|^2 \le \big(|\beta_1''|+|\beta_2''|\big)\|\tilde{v}_t\|^2+C\big(&\|\tilde{u}_{xx}\|^2+\|\tilde{u}_x\|^2+\|\tilde{v}_{xx}\|^2+\|\tilde{v}_x\|^2 +|\alpha'|\notag\\
&  +|\beta_1-\beta_2|+|\beta_1'|+|\beta_2'|+|\beta_1''|+|\beta_2''|\big).
\end{align}
Applying Gr\"onwall's inequality to \eqref{H24}--\eqref{H25}, using the assumptions of Theorem \ref{thm1} and Lemmas \ref{lem1}--\ref{lem3}, we can show that 
\begin{equation}\label{H26}
\|\tilde{u}_{t}(t)\|^2+\|\tilde{v}_{t}(t)\|^2+\int_0^t\left(\|\tilde{u}_{xt}(\tau)\|^2+ \varepsilon \|\tilde{v}_{xt}(\tau)\|^2\right)\mathrm{d}\tau \le C.
\end{equation}
As a consequence of the above estimate, we can show by using the equations in \eqref{e16} that 
\begin{equation*}
\|\tilde{u}_{xx}(t)\|^2+\|\tilde{v}_{xx}(t)\|^2+\int_0^t\left(\|\tilde{u}_{xxx}(\tau)\|^2+\varepsilon \|\tilde{v}_{xxx}(\tau)\|^2\right)\mathrm{d}\tau \le C.
\end{equation*}
We omit the routine technical details for brevity. This completes the proof of Lemma \ref{lem4}.
\end{proof}

Lemma \ref{lem1}--Lemma \ref{lem4} established the desired {\it a priori} estimates of the local solution. The global well-posedness of the IBVP then follows from these estimates and standard continuation argument. To finish the proof of Theorem \ref{thm1}, it remains to derive the time decay of the perturbation, which is carried out in the next subsection.

\subsection{Decay Estimate} 

\begin{lemma}\label{lem5}
Under the assumptions of Theorem $\ref{thm1}$, $\|\tilde{u}(t)\|_{H^2}+\|\tilde{v}(t)\|_{H^2} \to 0$, as $t\to\infty$. 
\end{lemma}

\begin{proof} {\bf Step 1.} From Lemma \ref{lem1}--Lemma \ref{lem2} we know that
\begin{align}\label{D0}
\|\tilde{u}_x(t)\|^2+\|\tilde{v}_x(t)\|^2 \in L^1(\mathbb{R}_+),
\end{align}
which, together with Poincar\'e's inequality, implies
\begin{align}\label{D1}
\|\tilde{u}(t)\|^2+\|\tilde{v}(t)\|^2 \in L^1(\mathbb{R}_+).
\end{align}
Since $0<\alpha(t)$ is uniformly bounded, we obtain from \eqref{D1} that
\begin{align}\label{D2}
\|\tilde{u}(t)\|^2+\alpha(t)\|\tilde{v}(t)\|^2 \in L^1(\mathbb{R}_+).
\end{align}
Using the assumptions of Theorem \ref{thm1} and Lemma \ref{lem1}--Lemma \ref{lem4}, we deduce from \eqref{e20} that
$$
\Big|\frac{\mathrm{d}}{\mathrm{d}t} \left(\|\tilde{u}\|^2+\alpha\|\tilde{v}\|^2\right) \Big| \le C\big(\|\tilde{u}_x\|^2+\|\tilde{v}_x\|^2 + |\beta_1-\beta_2| + |\alpha'| + |\beta_1'| + |\beta_2'|\big),
$$
where the constant is independent of $t$. Integrating the above inequality gives us
\begin{align}\label{D3}
\frac{\mathrm{d}}{\mathrm{d}t} \left(\|\tilde{u}(t)\|^2+\alpha(t)\|\tilde{v}(t)\|^2\right) \in L^1(\mathbb{R}_+).
\end{align}
From \eqref{D2} and \eqref{D3}, $\|\tilde{u}(t)\|^2+\alpha(t)\|\tilde{v}(t)\|^2 \in W^{1,1}(\mathbb{R}_+)$. Hence, $\|\tilde{u}(t)\|^2+\alpha(t)\|\tilde{v}(t)\|^2 \to 0$, as $t\to\infty$. Since $0<\underline{\alpha}\le \alpha(t)$, we conclude $\|\tilde{u}(t)\|^2+\|\tilde{v}(t)\|^2 \to 0$, as $t\to\infty$.

{\bf Step 2.} The decay of $\|\tilde{u}_x(t)\|^2+\|\tilde{v}_x(t)\|^2$ follows by the same idea. Indeed, from \eqref{e26} we have
\begin{align*}
&\,\Big|\frac{\mathrm{d}}{\mathrm{d}t} \left(\|\tilde{u}_x\|^2+\|\tilde{v}_x\|^2\right) \Big| \notag\\
\le&\, C\big(\|\tilde{u}_x\|^2+\|\tilde{v}_x\|^2 +\|\tilde{u}_{xx}\|^2+\|\tilde{v}_{xx}\|^2+ |\beta_1-\beta_2| + |\alpha'| + |\beta_1'| + |\beta_2'|\big),
\end{align*}
which, together with the estimates of the solution and the assumptions of Theorem \ref{thm1}, implies 
$$
\frac{\mathrm{d}}{\mathrm{d}t} \left(\|\tilde{u}_x\|^2+\|\tilde{v}_x\|^2\right) \in L^1(\mathbb{R}_+).
$$
Hence, by \eqref{D0}, we have $\|\tilde{u}_x(t)\|^2+\|\tilde{v}_x(t)\|^2 \in W^{1,1}(\mathbb{R}_+)$, which implies $\|\tilde{u}_x(t)\|^2+\|\tilde{v}_x(t)\|^2 \to 0$, as $t\to\infty$.

{\bf Step 3.} For the second order spatial derivatives, we know from \eqref{H26} and Poincar\'e's inequality that  
$$
\|\tilde{u}_t(t)\|^2+\|\tilde{v}_t(t)\|^2 \in L^1(\mathbb{R}_+).
$$
By \eqref{H21} and \eqref{H22a}, we can show that
$$
\begin{aligned}
\Big| \frac{\mathrm{d}}{\mathrm{d}t} \left(\|\tilde{u}_t\|^2+\|\tilde{v}_t\|^2\right) \Big| \le C\big(&\|\tilde{u}_t\|^2+\|\tilde{v}_t\|^2+\|\tilde{u}_{xt}\|^2+\|\tilde{v}_{xt}\|^2+\|\tilde{u}_x\|^2+\|\tilde{v}_x\|^2 \\
&+ |\alpha'| + |\alpha''| + |\beta_1-\beta_2| + |\beta_1'| + |\beta_2'| ++ |\beta_1''| + |\beta_2''|  \big),
\end{aligned}
$$
which, together with the estimates of the solution and the assumptions of Theorem \ref{thm1}, implies 
$$
\frac{\mathrm{d}}{\mathrm{d}t} \left(\|\tilde{u}_t(t)\|^2+\|\tilde{v}_t(t)\|^2\right) \in L^1(\mathbb{R}_+).
$$
Hence, $\|\tilde{u}_t(t)\|^2+\|\tilde{v}_t(t)\|^2 \in W^{1,1}(\mathbb{R}_+)$, which implies $\|\tilde{u}_t(t)\|^2+\|\tilde{v}_t(t)\|^2 \to 0$, as $t\to\infty$. According to \eqref{e16}, we can show that 
\begin{align}\label{D4}
&\,\|\tilde{u}_{xx}\|^2+\|\tilde{v}_{xx}\|^2 \notag\\
\le &\,C\big( \|\tilde{u}_t\|^2 + \|\tilde{v}_t\|^2 + \|\tilde{u}_x\|^2 + \|\tilde{v}_x\|^2 + |\alpha'| + |\beta_1-\beta_2| + |\beta_1'| + |\beta_2'|\big). 
\end{align}
Since $\alpha'$, $\beta_1-\beta_2$, $\beta_1'$ and $\beta_2'$ are assumed to belong to $W^{1,1}(\mathbb{R}_+)$, they all tend to zero as $t\to\infty$. Therefore, the decay of $\|\tilde{u}_{xx}(t)\|^2+\|\tilde{v}_{xx}(t)\|^2$ follows from \eqref{D4} and the decay of the first order derivatives of the perturbation. This completes the proof of Lemma \ref{lem5}. 
\end{proof}

\section{Proof of Theorem \ref{thm2}}\label{PoT2}

This section is devoted to prove Theorem \ref{thm2}. Recall  system \eqref{1.2}:
\begin{subequations}\label{3.1}
\begin{alignat}{2}
u_t-(uv)_x &=u_{xx}, \label{3.1a}\\
v_t-u_x&= 0, \label{3.1b} 
\end{alignat}
\end{subequations}
together with the initial and boundary conditions:
\begin{subequations}\label{3.1a}
\begin{alignat}{2}
(u,v)(x,0) &= (u_0,v_0)(x), \quad&x&\in [0,1],\label{3.1aa}\\
u(0,t)&=u(1,t)=\alpha(t),\quad &t& \ge0. \label{3.1ab}
\end{alignat}
\end{subequations}
We stress that some of the {\it a priori} estimates in the previous section cannot be carried over to \eqref{3.1}. This can be seen from \eqref{e27}, in which the constant blows up when $\varepsilon\to0$. We adapt a different approach to estimate higher order spatial derivatives of the solution.

{\bf Step 1.} First we note that \eqref{e5} is still valid when $\varepsilon=0$. We record it here for convenience:
\begin{equation}\label{3.0a}
\begin{aligned}
\frac{\mathrm{d}}{\mathrm{d}t}E(u,\alpha)+\int_0^1\frac{(u_x)^2}{u}\mathrm{d}x = -\int_0^1 v u_x\mathrm{d}x-\int_0^1(u-\alpha)\frac{\alpha'}{\alpha}\mathrm{d}x,
\end{aligned}
\end{equation}
where $E(u,\alpha)$ is the same as before. Integrating \eqref{3.1b} from $x=0$ to $x=1$ and using the boundary conditions for $u$, we have
\begin{align}\label{3.2}
\frac{\mathrm{d}}{\mathrm{d}t}\int_0^1 v(x,t)\dx = 0.
\end{align}
Integrating \eqref{3.2} with respect to $t$, we obtain
\begin{align}\label{3.3}
\int_0^1 v(x,t)\dx = \int_0^1 v(x,0)\dx \equiv \overline{v}.
\end{align}
Since $\overline{v}$ is a constant, we easily obtain
\begin{align}\label{3.0b}
(v-\overline{v})_t - u_x=0.
\end{align}
Taking $L^2$ inner product of \eqref{3.0b} with $(v-\overline{v})$, we have
\begin{align}\label{3.0c}
\frac12\frac{\mathrm{d}}{\mathrm{d}t} \|v-\overline{v}\|^2 = \int_0^1 u_x(v-\overline{v})\dx = \int_0^1 u_x v\dx,
\end{align}
where the integral of $u_x\overline{v}$ is zero, since $\overline{v}$ is a constant and the boundary values of $u$ match at the endpoints. Taking the sum of \eqref{3.0a} and \eqref{3.0c} gives us
\begin{equation}\label{3.0d}
\begin{aligned}
\frac{\mathrm{d}}{\mathrm{d}t} \Big(E(u,\alpha) + \frac12\|v-\overline{v}\|^2\Big) +\int_0^1\frac{(u_x)^2}{u}\mathrm{d}x = -\int_0^1(u-\alpha)\frac{\alpha'}{\alpha}\mathrm{d}x.
\end{aligned}
\end{equation}
Similar to {\bf Step 2} of the proof of Lemma \ref{lem1}, we can show that
\begin{equation}\label{3.0e}
\begin{aligned}
\frac{\mathrm{d}}{\mathrm{d}t} \Big(E(u,\alpha) + \frac12\|v-\overline{v}\|^2\Big) +\int_0^1\frac{(u_x)^2}{u}\mathrm{d}x \le \frac{|\alpha'|}{\underline{\alpha}}E(u,\alpha) + e|\alpha'|.
\end{aligned}
\end{equation}
Applying Gr\"onwall's inequality and using the assumptions for $\alpha$, we obtain
\begin{align}\label{3.0f}
E(u,\alpha)(t)+\|(v-\overline{v})(t)\|^2+\int_0^t\int_0^1\frac{(u_x)^2}{u}\mathrm{d}x\mathrm{d}\tau \le C,
\end{align}
where the constant on the right-hand side is independent of $t$.

{\bf Step 2.} Again, we define the perturbed variables:
$$
\tu(x,t) = u(x,t) - \alpha(t),\quad \tv(x,t) = v(x,t) - \overline{v}.
$$
Similar to \eqref{e16}, we have
\begin{subequations}\label{3.4}
\begin{alignat}{2}
\tilde{u}_t-(\tilde{u}\tilde{v})_x &=\tilde{u}_{xx} + \alpha \tilde{v}_x + \overline{v}\tilde{u}_x - \alpha', \label{3.4a}\\
\tilde{v}_t-\tilde{u}_x &= 0.\label{3.4b} 
\end{alignat}
\end{subequations}
Taking $L^2$ inner product of \eqref{3.4a} with $\tu$, we have
\begin{align}\label{3.5}
\frac12\frac{\mathrm{d}}{\mathrm{d}t} \|\tu\|^2 + \|\tu_x\|^2 = -\int_0^1 \tu\tv\tu_x\dx + \alpha\int_0^1 \tv_x\tu\dx - \alpha'\int_0^1 \tu\dx.
\end{align}
Taking $L^2$ inner product of \eqref{3.4b} with $\alpha \tv$, we obtain
\begin{align}\label{3.6}
\frac12\frac{\mathrm{d}}{\mathrm{d}t}\big(\alpha\|\tv\|^2\big) = \alpha\int_0^1 \tu_x\tv\dx + \frac{\alpha'}{2}\|\tv\|^2.
\end{align}
Taking the sum of \eqref{3.5} and \eqref{3.6} gives us
\begin{align}\label{3.7}
\frac12\frac{\mathrm{d}}{\mathrm{d}t} \big(\|\tu\|^2 +\alpha\|\tv\|^2\big) + \|\tu_x\|^2 = -\int_0^1 \tu\tv\tu_x\dx - \alpha'\int_0^1 \tu\dx + \frac{\alpha'}{2}\|\tv\|^2.
\end{align}
The estimate of the first integral on the right-hand side of \eqref{3.7} is identical to \eqref{e23a}, which is recorded here for convenience:
\begin{equation}\label{3.7a}
\int_0^1 \tu\tv\tu_x\dx \le C\int_0^1\frac{(u_x)^2}{u}\mathrm{d}x+\frac12\|\tilde{u}_x\|^2.
\end{equation}
The second integral on the right-hand side of \eqref{3.7} is estimated as 
\begin{align}\label{3.7b}
\alpha'\int_0^1 \tu\dx \le |\alpha'|\,\|\tu\| \le \frac{|\alpha'|}{2}\|\tu\|^2 + \frac{|\alpha'|}{2}.
\end{align}
Substituting \eqref{3.7a} and \eqref{3.7b} into \eqref{3.7}, we obtain
\begin{align}\label{3.7c}
\frac12\frac{\mathrm{d}}{\mathrm{d}t} \big(\|\tu\|^2 +\alpha\|\tv\|^2\big) + \frac12 \|\tu_x\|^2 \le C|\alpha'|\big(\|\tu\|^2 +\alpha\|\tv\|^2\big) + C\Big(\int_0^1\frac{(u_x)^2}{u}\mathrm{d}x + |\alpha'|\Big),
\end{align}
where we utilized the strictly positive lower bound of $\alpha$ for the third term on the right-hand side of \eqref{3.7}. Applying Gr\"onwall's inequality and using the assumptions for $\alpha$, we can show that
\begin{align}\label{3.7d}
\|\tu(t)\|^2 + \|\tv(t)\|^2 + \int_0^t\|\tu_x(\tau)\|^2\mathrm{d}\tau \le C,
\end{align}
where \eqref{3.0f} is applied.

{\bf Step 3.} This step is the major difference between the proof of Theorem \ref{thm1} and Theorem \ref{thm2}. The following equation, obtained by substituting $\tu_{xx}=\tv_{xt}$ into \eqref{3.4a}, 
\begin{align}\label{3.8}
\tv_{xt} + \alpha \tilde{v}_x = \tilde{u}_t-(\tilde{u}\tilde{v})_x - \overline{v}\tilde{u}_x + \alpha',
\end{align}
serves as the foundation for obtaining the time integrability of the Sobolev norm of $\tv$. Taking $L^2$ inner product of \eqref{3.8} with $\tv_x$, we have
\begin{align}\label{3.9}
\frac12\frac{\mathrm{d}}{\mathrm{d}t} \|\tv_x\|^2 + \alpha\|\tv_x\|^2 = \int_0^1 \tu_t\tv_x\dx - \int_0^1 (\tu\tv)_x\tv_x\dx - \overline{v}\int_0^1 \tu_x\tv_x\dx + \alpha'\int_0^1\tv_x\dx.
\end{align}
Using \eqref{3.4b}, we rewrite the first integral on the right-hand side of \eqref{3.9} as 
\begin{align}\label{3.10}
\int_0^1 \tu_t\tv_x\dx = \frac{\mathrm{d}}{\mathrm{d}t}\int_0^1 \tu\tv_x\dx - \int_0^1 \tu\tv_{xt}\dx &= \frac{\mathrm{d}}{\mathrm{d}t}\int_0^1 \tu\tv_x\dx - \int_0^1 \tu\tu_{xx}\dx \notag\\
&= \frac{\mathrm{d}}{\mathrm{d}t}\int_0^1 \tu\tv_x\dx + \|\tu_x\|^2.
\end{align}
Substituting \eqref{3.10} into \eqref{3.9}, we obtain
\begin{align}\label{3.11}
&\frac{\mathrm{d}}{\mathrm{d}t} \Big(\frac12\|\tv_x\|^2 -\int_0^1\tu\tv_x\dx\Big) + \alpha\|\tv_x\|^2\notag\\
 = &\|\tu_x\|^2 - \int_0^1 (\tu\tv)_x\tv_x\dx - \overline{v}\int_0^1 \tu_x\tv_x\dx + \alpha'\int_0^1\tv_x\dx.
\end{align}
The three integrals on the right-hand side of \eqref{3.11} are estimated as
\begin{align*}
\bullet\quad&\int_0^1 (\tu\tv)_x\tv_x\dx \le C\|\tu_x\|\|\tv_x\|\tv_x\| \le \frac{\alpha}{6}\|\tv_x\|^2 + \frac{C}{\alpha}\|\tu_x\|^2\|\tv_x\|^2, \notag\\
\bullet\quad&\overline{v}\int_0^1 \tu_x\tv_x\dx \le C\|\tu_x\| \|\tv_x\| \le \frac{\alpha}{6}\|\tv_x\|^2 + \frac{C}{\alpha}\|\tu_x\|^2, \notag\\
\bullet\quad&\alpha'\int_0^1\tv_x\dx \le |\alpha'| \|\tv_x\| \le \frac{\alpha}{6}\|\tv_x\|^2 + \frac{C}{\alpha}|\alpha'|^2.
\end{align*}
Substituting these estimates into \eqref{3.11} gives us
\begin{align}\label{3.12}
\frac{\mathrm{d}}{\mathrm{d}t} \Big(\frac12\|\tv_x\|^2 -\int_0^1\tu\tv_x\dx\Big) + \frac{\alpha}{2}\|\tv_x\|^2 \le C\|\tu_x\|^2\|\tv_x\|^2 + C\big( \|\tu_x\|^2 +|\alpha'|\big),
\end{align}
where we used the information that $\alpha(t)\ge \underline{\alpha}>0$ and $|\alpha'|$ is uniformly bounded. According to \eqref{3.7d}, there is a constant, denoted by $\overline{U}_1$, such that $\|\tu\|^2 \le \overline{U}_1$. Since 
\begin{align*}
\Big| \int_0^1 \tu\tv_x\dx \Big| &\le \|\tu\|^2 + \frac14\|\tv_x\|^2 \le \overline{U}_1 + \frac14\|\tv_x\|^2, 
\end{align*}
we know that
\begin{align}\label{3.13}
E_1(t) := \frac12\|\tv_x\|^2 -\int_0^1\tu\tv_x\dx + \overline{U}_1 \ge \frac14\|\tv_x\|^2.
\end{align}
Then we update \eqref{3.12} as 
\begin{align}\label{3.14}
\frac{\mathrm{d}}{\mathrm{d}t} E_1(t) + \frac{\alpha}{2}\|\tv_x\|^2 \le C\|\tu_x\|^2 E_1(t) + C\big( \|\tu_x\|^2 +|\alpha'|\big).
\end{align}
Applying Gr\"onwall's inequality,  using \eqref{3.7d} and the assumptions for $\alpha$, we can show that
\begin{align}\label{3.15}
\|\tv_x(t)\|^2 + \int_0^t \|\tv_x(\tau)\|^2\mathrm{d}\tau \le C,
\end{align}
where \eqref{3.13} is applied and the constant on the right-hand side is independent of $t$.

{\bf Step 4.} Taking $L^2$ inner product of \eqref{3.4a} with $-\tu_{xx}$, we have 
\begin{align}\label{3.16}
&\frac12\frac{\mathrm{d}}{\mathrm{d}t} \|\tu_x\|^2 + \|\tu_{xx}\|^2\notag\\
 = & -\int_0^1 (\tu\tv)_x\tu_{xx}\dx - \alpha\int_0^1 \tv_x\tu_{xx}\dx - \overline{v}\int_0^1 \tu_x\tu_{xx}\dx + \alpha'\int_0^1 \tu_{xx}\dx.
\end{align}
The integrals on the right-hand side of \eqref{3.16} are estimated as
\begin{align*}
\bullet\quad&\int_0^1 (\tu\tv)_x\tu_{xx}\dx \le C\|\tu_x\|\|\tv_x\| \|\tu_{xx}\|,\\
\bullet\quad&\alpha\int_0^1 \tv_x\tu_{xx}\dx \le C\|\tv_x\| \|\tu_{xx}\|,\\
\bullet\quad&\overline{v}\int_0^1 \tu_x\tu_{xx}\dx \le \overline{v}\|\tu_x\| \|\tu_{xx}\|,\\
\bullet\quad&\alpha'\int_0^1 \tu_{xx}\dx \le |\alpha'| \|\tu_{xx}\|.
\end{align*}
Substituting these estimates into \eqref{3.16} and applying the Cauchy-Schwarz inequality, we obtain
\begin{align}\label{3.17}
\frac12\frac{\mathrm{d}}{\mathrm{d}t} \|\tu_x\|^2 + \frac12\|\tu_{xx}\|^2 \le C\|\tu_x\|^2\|\tv_x\|^2 + C\big(\|\tu_x\|^2 + \|\tv_x\|^2 + |\alpha'|\big),
\end{align}
where the uniform boundedness of $|\alpha'|$ is used. Applying Gr\"onwall's inequality to \eqref{3.17}, using \eqref{3.7d}, \eqref{3.15} and the assumptions for $\alpha$, we can show that
\begin{align}\label{3.18}
\|\tu_x(t)\|^2 + \int_0^t\|\tu_{xx}(\tau)\|^2\mathrm{d}\tau \le C,
\end{align}
where the constant on the right-hand side is independent of $t$.

{\bf Step 5.} Differentiating \eqref{3.4a} with respect to $t$, we have
\begin{align}\label{3.19}
\tilde{u}_{tt}-(\tilde{u}\tilde{v})_{xt}=\tilde{u}_{xxt} + (\alpha \tilde{v})_{xt} + \overline{v}\tilde{u}_{xt} - \alpha''.
\end{align}
Taking $L^2$ inner product of \eqref{3.19} with $\tu_t$ and integrating by parts, we obtain
\begin{align}\label{3.20}
\frac12\frac{\mathrm{d}}{\mathrm{d}t}\|\tilde{u}_t\|^2 + \|\tu_{xt}\|^2 = -\int_0^1 (\tilde{u}\tilde{v})_{t}\tu_{xt}\dx - \int_0^1 (\alpha \tilde{v})_{t} \tu_{xt}\dx - \alpha'' \int_0^1 \tu_t \dx.
\end{align}
The integrals on the right-hand side of \eqref{3.20} are estimated as:
\begin{align*}
\bullet\quad&\int_0^1 (\tilde{u}\tilde{v})_{t}\tu_{xt}\dx \le C\big(\|\tu_x\|\|\tv_t\| + \|\tv_x\| \|\tu_t\|\big)\|\tu_{xt}\|, \\
\bullet\quad&\int_0^1 (\alpha \tilde{v})_{t} \tu_{xt}\dx \le \big(|\alpha'|\|\tv\| + |\alpha| \|\tv_t\|\big) \|\tu_{xt}\|, \\
\bullet\quad&\alpha'' \int_0^1 \tu_t \dx \le \frac{|\alpha''|}{2} +  \frac{|\alpha''|}{2}\|\tu_t\|^2.
\end{align*}
Using these estimates and the Cauchy-Schwarz inequality, we update \eqref{3.20} as 
\begin{align}\label{3.21}
&\,\frac12\frac{\mathrm{d}}{\mathrm{d}t}\|\tilde{u}_t\|^2 + \frac12\|\tu_{xt}\|^2 \notag\\
\le &\,C\big(\|\tv_x\|^2+|\alpha''|\big) \|\tu_t\|^2 + C\big(\|\tu_x\|^2\|\tv_t\|^2 + \|\tv\|^2 + \|\tv_t\|^2 + |\alpha''|\big) \notag \\
\le &\,C\big(\|\tv_x\|^2+|\alpha''|\big) \|\tu_t\|^2 + C\big(\|\tu_x\|^2 + \|\tv_x\|^2 + |\alpha''|\big),
\end{align}
where we replaced $\tu_x$ by $\tv_t$, utilized \eqref{3.18} and the uniform boundedness of $|\alpha'|$ and $|\alpha|$, together with Poincar\'e's inequality for $\tv$ (since it is mean free). Applying Gr\"onwall's inequality gives us
\begin{align}\label{3.22}
\|\tilde{u}_t(t)\|^2 + \int_0^t \|\tu_{xt}(\tau)\|^2\mathrm{d}\tau \le C, 
\end{align}
where the constant on the right-hand side is independent of $t$. As a consequence of \eqref{3.22}, \eqref{3.4a} and previous estimates, we have 
\begin{align}\label{3.22a}
\|\tu_{xx}(t)\|^2 \le C\big(\|\tu_t\|^2 + \|\tu_x\|^2\|\tv_x\|^2 + \|\tv_x\|^2 + \|\tu_x\|^2 + |\alpha'|^2\big) \le C,
\end{align}
where the constant $C$ is independent of $t$.

{\bf Step 6.} Taking $\partial_x$ of \eqref{3.8} gives us
\begin{align}\label{3.23}
\tv_{xxt} + \alpha \tilde{v}_{xx} = \tilde{u}_{xt}-(\tilde{u}\tilde{v})_{xx} - \overline{v}\tilde{u}_{xx}.
\end{align}
Taking $L^2$ inner product of \eqref{3.23} with $\tv_{xx}$, we have
\begin{align}\label{3.24}
\frac12\frac{\mathrm{d}}{\mathrm{d}t}\|\tv_{xx}\|^2 + \alpha \|\tilde{v}_{xx}\| = \int_0^1 \tilde{u}_{xt} \tv_{xx}\dx -\int_0^1 (\tilde{u}\tilde{v})_{xx} \tv_{xx}\dx - \overline{v} \int_0^1 \tilde{u}_{xx} \tv_{xx}\dx.
\end{align}
The integrals on the right-hand side are estimated as:
\begin{align*}
\bullet\quad&\int_0^1 \tilde{u}_{xt} \tv_{xx}\dx \le C\|\tu_{xt}\|^2 + \frac{\alpha}{6}\|\tv_{xx}\|^2, \\
\bullet\quad&\int_0^1 (\tilde{u}\tilde{v})_{xx} \tv_{xx}\dx \le C\big(\|\tv_x\|^2\|\tu_{xx}\|^2 + \|\tu_x\|^2\|\tv_{xx}\|^2\big) + \frac{\alpha}{6}\|\tv_{xx}\|^2, \\
\bullet\quad&\overline{v} \int_0^1 \tilde{u}_{xx} \tv_{xx}\dx \le C\|\tu_{xx}\|^2 + \frac{\alpha}{6}\|\tv_{xx}\|^2.
\end{align*}
Then we update \eqref{3.24} as 
\begin{align}\label{3.25}
\frac12\frac{\mathrm{d}}{\mathrm{d}t}\|\tv_{xx}\|^2 + \frac{\alpha}{2} \|\tilde{v}_{xx}\| \le C\|\tu_x\|^2\|\tv_{xx}\|^2 + C\big(\|\tu_{xt}\|^2 + \|\tu_{xx}\|^2\big),
\end{align}
where we used \eqref{3.15} for the estimate of $\|\tv(t)\|$. Applying Gr\"onwall's inequality, we obtain
\begin{align}\label{3.26}
\|\tv_{xx}(t)\|^2 + \int_0^t \|\tv_{xx}(\tau)\|^2\mathrm{d}\tau \le C, 
\end{align}
where the constant on the right-hand side is independent of $t$. As a consequence of \eqref{3.26}, \eqref{3.4a} and previous estimates, we can show that
\begin{align}\label{3.27}
\int_0^t \|\tu_{xxx}(\tau)\|^2\mathrm{d}\tau \le C,
\end{align}
where the constant on the right-hand side is independent of $t$. The routine technical details are omitted for brevity. This completes the proof of the {\it a priori} estimates of the solution, as stated in Theorem \ref{thm2}. Moreover, the time decay of the perturbation follows from the same spirit in the previous section. The proof of Theorem \ref{thm2} is thus complete.

\section{Numerical Simulations}\label{NS}
In this section, we carry out numerical tests to study the dynamical properties of solutions to \eqref{1.1}, and those of \eqref{1.2}. In particular, we provide numerical confirmation of the rigorous qualitative results established above (cf. Theorem \ref{thm1} and Theorem \ref{thm2}), and explore the steady state of the solutions when relaxing some of the assumptions of the above mentioned theorems. 

We perform numerical simulations to the transformed systems \eqref{1.1} and  \eqref{1.2} given the singularity of the term $(\log c)_x$ in the original system.  We make use of an explicit finite-difference scheme to solve the equations with a second order approximation of spatial derivatives with a temporal mesh, $\Delta t= \frac{(\Delta x)^2}2$, and a spatial mesh, $\Delta x<\sqrt{\frac{\varepsilon}{10}}$, similar to what would be chosen for the heat equation in order to avoid that the numerical diffusion dominates the chemical diffusion. The domain for the system is the unit interval $[0,1]$ with time--dependant boundary data $\alpha_1, \alpha_2$ and $\beta_1, \beta_2$. For the following results, we have used initial data $u(x,0)=0.3+0.1 \sin^2(2\pi x)$ and \padi{$v(x,0)=0.2 \sin^2(2\pi x)$.} We would like to point out that the large--time behaviors of the solutions do not vary for similar kinds of initial data. For consistency, we use the same initial conditions and $N=200$ throughout the computations in this section. 

We define $A (x,t)= (\alpha_2-\alpha_1)x+ \alpha_1$ and $B (x,t)= (\beta_2-\beta_1)x+ \beta_1$, the linear interpolation of the boundary data. Notice that $\tilde u=u-A$ and $\tilde v=v-B$. Throughout this section, the letters with overhead bar denote real constants. 

\subsection{Study of solutions of \eqref{1.1}} We study the existence of  steady states of  solutions  of \eqref{1.1} for different conditions on the boundary data. For these computations, we use $\varepsilon=0.7$. 
\begin{itemize}
\item \textbf{Case 1}  $\boxed{\alpha_1=\alpha_2\to\bar{\alpha}, \ \  \beta_1\neq\beta_2,\ \ \beta_1\to\bar\beta\gets\beta_2}$ 
\end{itemize}
The solution  $(u,v)$  converges to the steady state $(\bar\alpha,\bar\beta)$ as proved in  Theorem \ref{thm1}.\\
\begin{figure}[htbp]
    \begin{subfigure}{0.45\textwidth}
        \includegraphics[width=\textwidth]{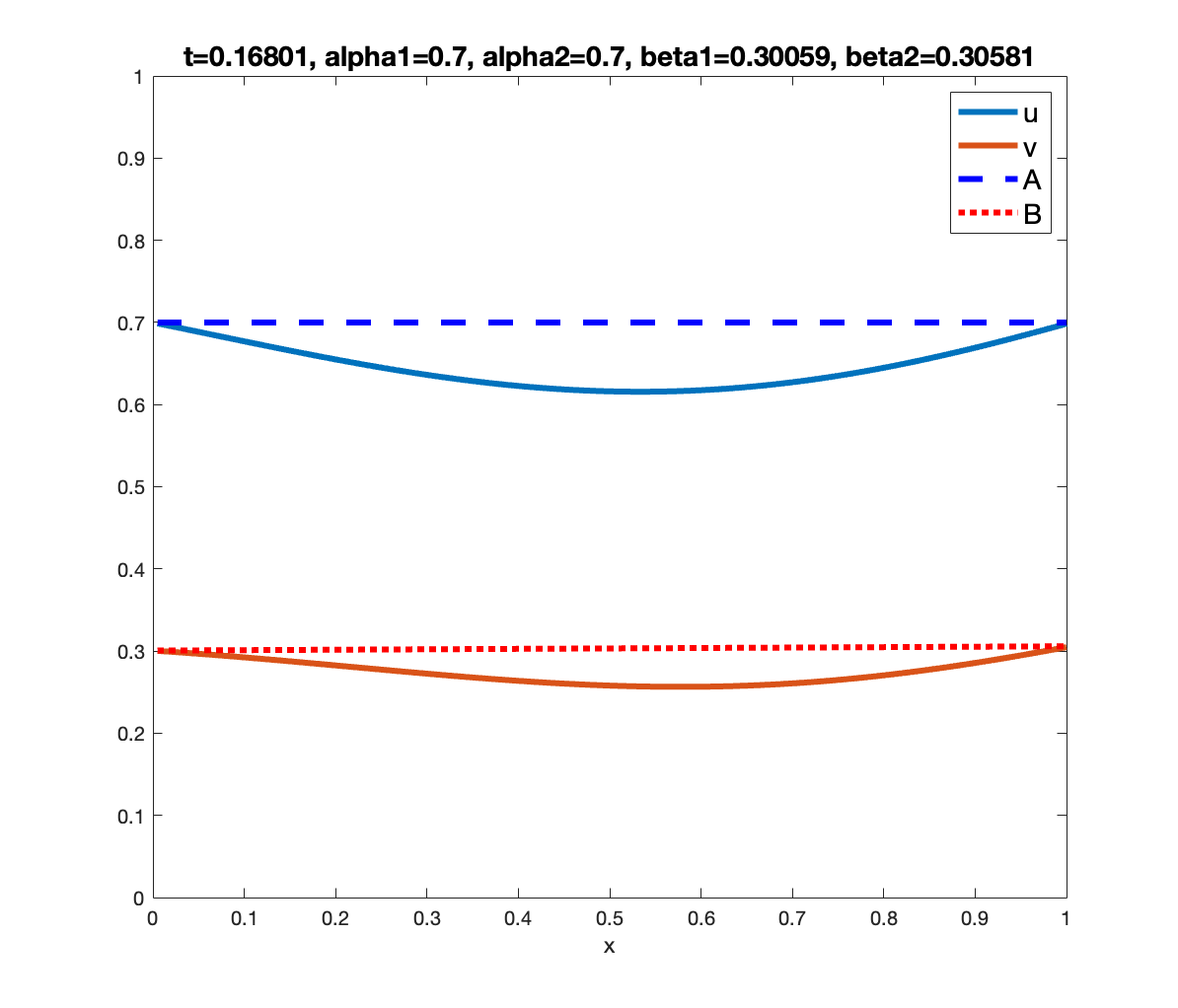}
        \caption{$t=0.16801$}
        \label{fig1}
    \end{subfigure}
    \qquad
    \begin{subfigure}{0.45\textwidth}
        \includegraphics[width=\textwidth]{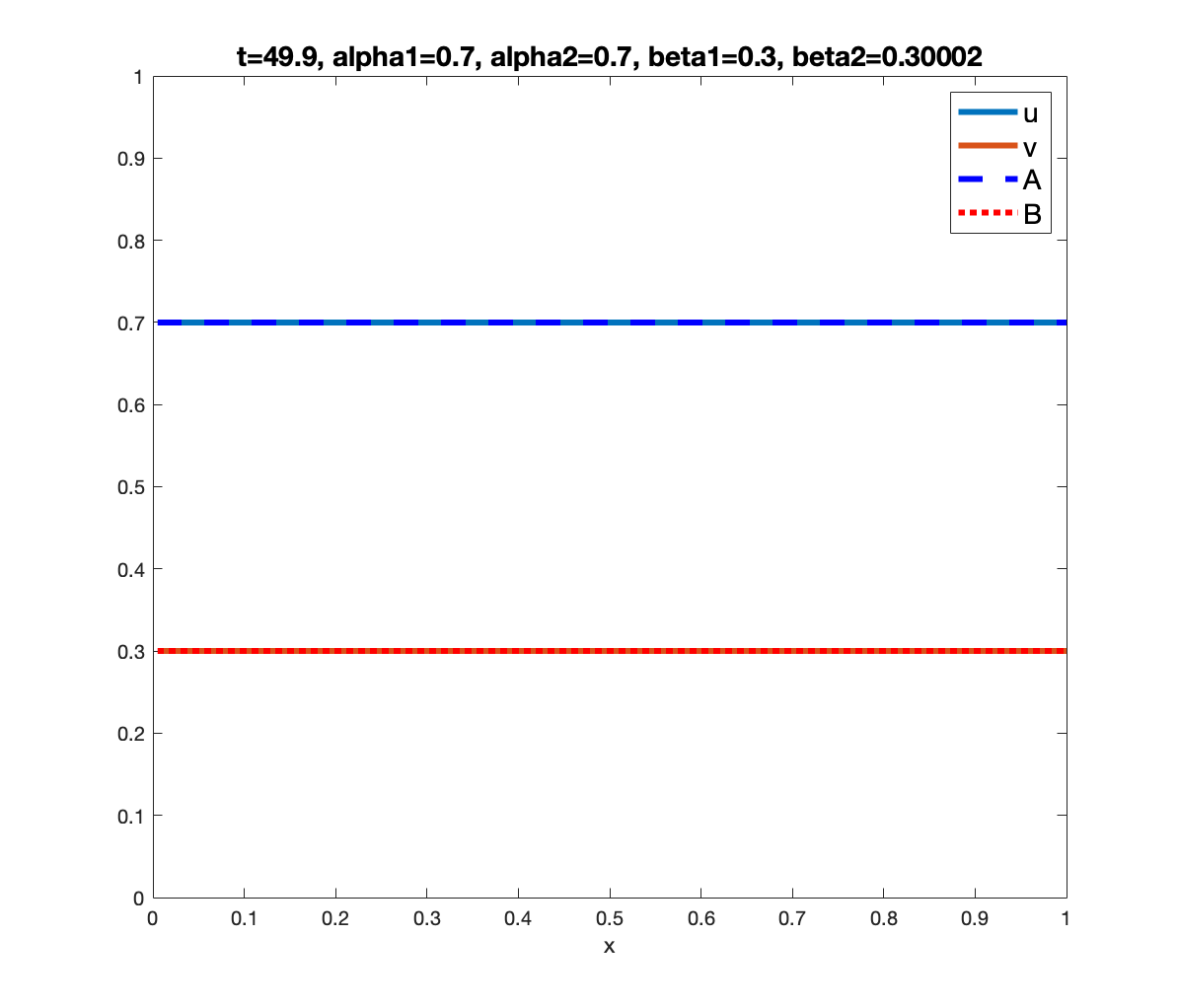}
        \caption{$t=49.9$}
        \label{fig2}
    \end{subfigure}
    \label{case1}
    \caption{Solution of \eqref{1.1} with $\alpha_1=\alpha_2\to\bar\alpha, \ \  \beta_1\neq\beta_2,\ \ \beta_1\to\bar\beta\gets\beta_2$.}
\end{figure}\\
In this case, we plot $u,v, A, B$ with $\alpha_1=0.7+\exp(-200000t)=\alpha_2$,  $\beta_1=0.3+1/(1+10000t)$ and $\beta_2= 0.3+1/(4+1000t)$.
Figure \ref{fig1} plots the solution $(u,v)$ and the linear interpolation $(A,B)$ at time $t=0.168$, with  $\alpha_1=0.7=\alpha_2$, $\beta_1=0.301, \beta_2=0.306$. Figure \ref{fig2} plots $(u,v,A,B)$ at time $t=49.9$, with  $\alpha_1=0.7=\alpha_2$, $\beta_1=0.3002, \beta_2=0.3018$. In particular, we observe that the solution reaches the steady state $(\bar\alpha,\bar\beta)=(0.7,0.3)$, as $t\to\infty$.

\begin{itemize}
\item \textbf{Case 2}  $\boxed{\alpha_1\neq\alpha_2, \ \ \alpha_1\to\bar\alpha\gets\alpha_2,\ \ \beta_1\neq\beta_2,\ \ \beta_1\to\bar\beta\gets\beta_2}$ 
\end{itemize}
The solution  $(u,v)$  converges to the steady state $(\bar\alpha,\bar\beta)$.\\
\begin{figure}[htbp]
    \begin{subfigure}{0.45\textwidth}
        \includegraphics[width=\textwidth]{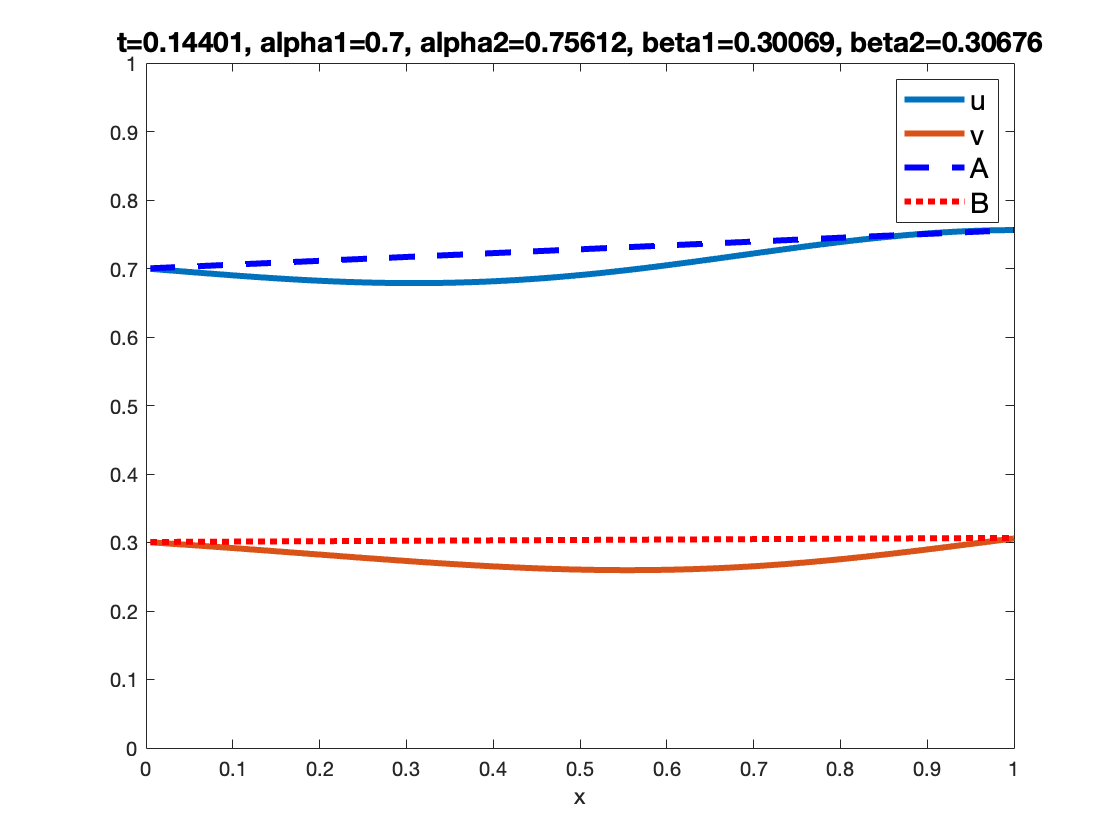}
        \caption{$t=0.1440$}
        \label{figex1}
    \end{subfigure}
    \qquad
    \begin{subfigure}{0.4\textwidth}
        \includegraphics[width=\textwidth]{fig2.png}
        \caption{$t=49.99$}
        \label{figex2}
    \end{subfigure}
    \label{casenew2}
    \caption{Solution of \eqref{1.1} with $\alpha_1\neq\alpha_2, \ \ \alpha_1\to\bar\alpha\gets\alpha_2,\ \ \beta_1\neq\beta_2,\ \ \beta_1\to\bar\beta\gets\beta_2$.}
\end{figure}\\
In this case, we plot $u,v, A, B$ with $\alpha_1=0.7+\exp(-200000t)$, $\alpha_2=0.7+\exp(-20t)$,  $\beta_1=0.3+1/(1+10000t)$ and $\beta_2= 0.3+1/(4+1000t)$.
Figure \ref{figex1} plots the solution $(u,v)$ and the linear interpolation $(A,B)$ at time $t=0.1440$, with  $\alpha_1=0.7$, $\alpha_2=0.756$, $\beta_1=0.300,\, \beta_2=0.306$. Figure \ref{figex2} plots $(u,v,A,B)$ at time $t=49.99$, with  $\alpha_1=0.7=\alpha_2$, $\beta_1=0.3, \beta_2=0.3002$. In particular, we observe that the solution reaches the steady state $(\bar\alpha,\bar\beta)=(0.7,0.3)$ as in \textbf{Case 1}.

\begin{itemize}
\item \textbf{Case 3}  $\boxed{\alpha_1=\alpha_2\to\bar\alpha, \ \  \beta_1\to\bar\beta_1\neq\bar\beta_2\gets\beta_2}$ 
\end{itemize}
The solution $(u,v)$  has a steady state different from $\big(\bar\alpha,(\bar\beta_2-\bar\beta_1)x+\bar\beta_1\big)$.\\
\begin{figure}[htbp]
    \begin{subfigure}{0.45\textwidth}
        \includegraphics[width=\textwidth]{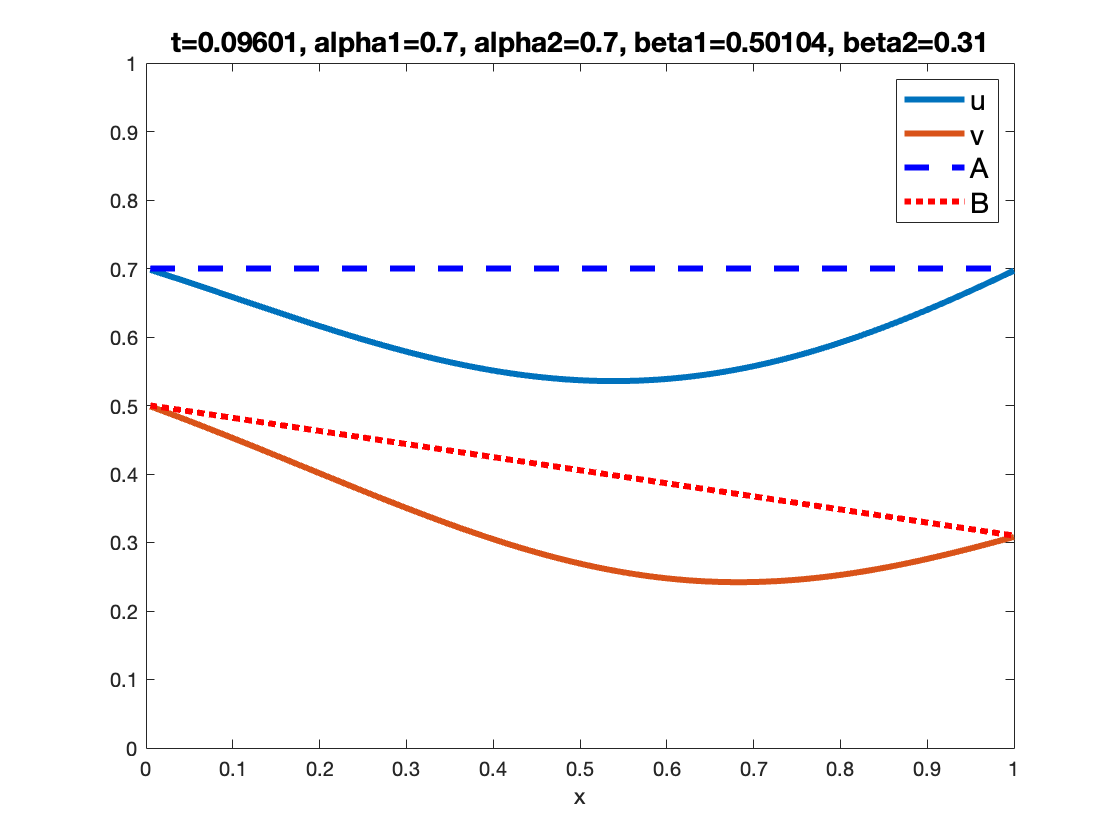}
        \caption{$t=0.09601$}
        \label{fig3}
    \end{subfigure}
    \qquad
    \begin{subfigure}{0.45\textwidth}
        \includegraphics[width=\textwidth]{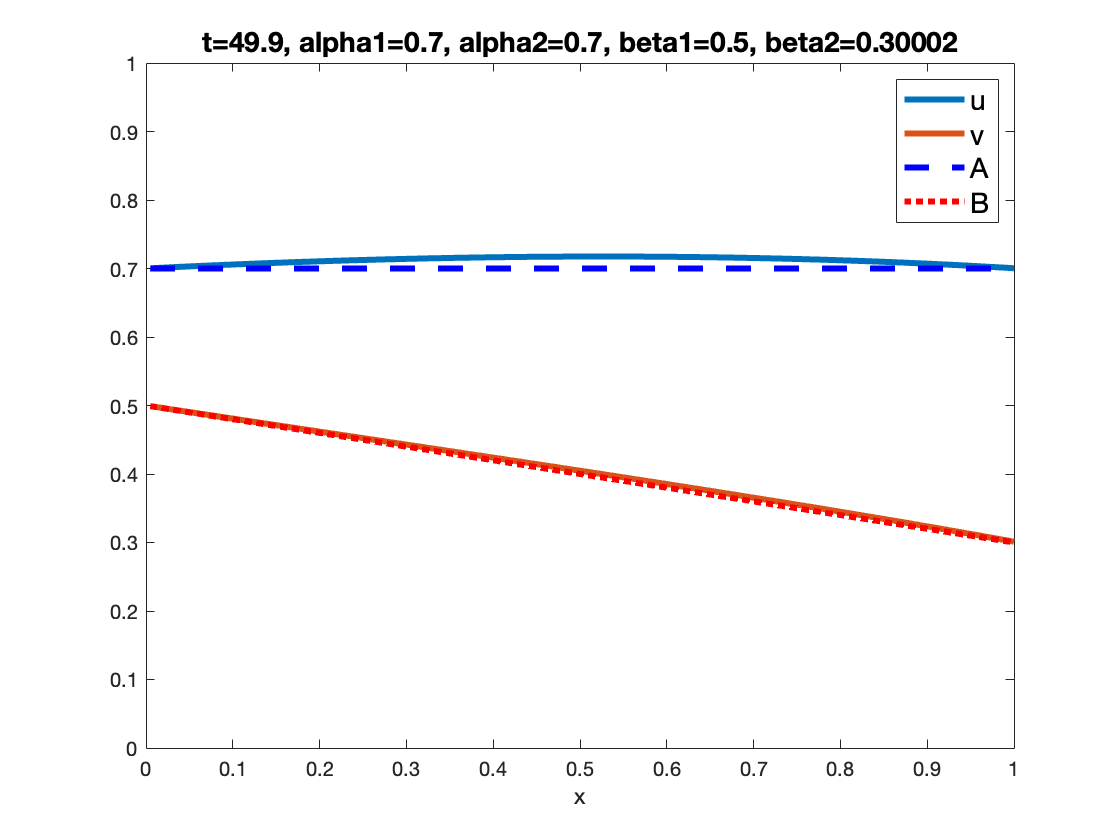}
        \caption{$t=49.9$}
        \label{fig4}
    \end{subfigure}
    \caption{Solution of \eqref{1.1} with $\alpha_1=\alpha_2\to\bar\alpha, \ \  \beta_1\to\bar\beta_1\neq\bar\beta_2\gets\beta_2.$}
    \label{case2}
\end{figure}\\
In this case, we plot $u,v, A, B$ with $\alpha_1=0.7+\exp(-200000t)=\alpha_2$,  $\beta_1=0.5+1/(1+10000t)$ and $\beta_2= 0.3+1/(4+1000t)$.
Figure \ref{fig3} plots the solution $(u,v)$ and the linear interpolation $(A,B)$ at time $t=0.168$, with  $\alpha_1=0.7=\alpha_2$, $\beta_1=0.301, \beta_2=0.306$. Figure \ref{fig4} shows $(u,v,A,B)$ at time $t=49.9$, with  $\alpha_1=0.7=\alpha_2$, $\beta_1=0.500, \beta_2=0.300$. In particular, we observe that the solution reaches a steady state different from $\big(\bar\alpha,(\bar\beta_2-\bar\beta_1)x+\bar\beta_1\big)=\big(0.7,-0.2x+0.5\big)$.
 \begin{itemize}
\item \textbf{Case 4}  $\boxed{\alpha_1\to\bar\alpha_1\neq\bar\alpha_2\gets\alpha_2,\ \ \beta_1\to\bar\beta_1\neq\bar\beta_2\gets\beta_2}$ 
\end{itemize}
The solution  $(u,v)$  has a steady state different from $\big((\bar\alpha_2-\bar\alpha_1)x+\bar\alpha_1,(\bar\beta_2-\bar\beta_1)x+\bar\beta_1\big)$.\\
\begin{figure}[htbp]
    \begin{subfigure}{0.45\textwidth}
        \includegraphics[width=\textwidth]{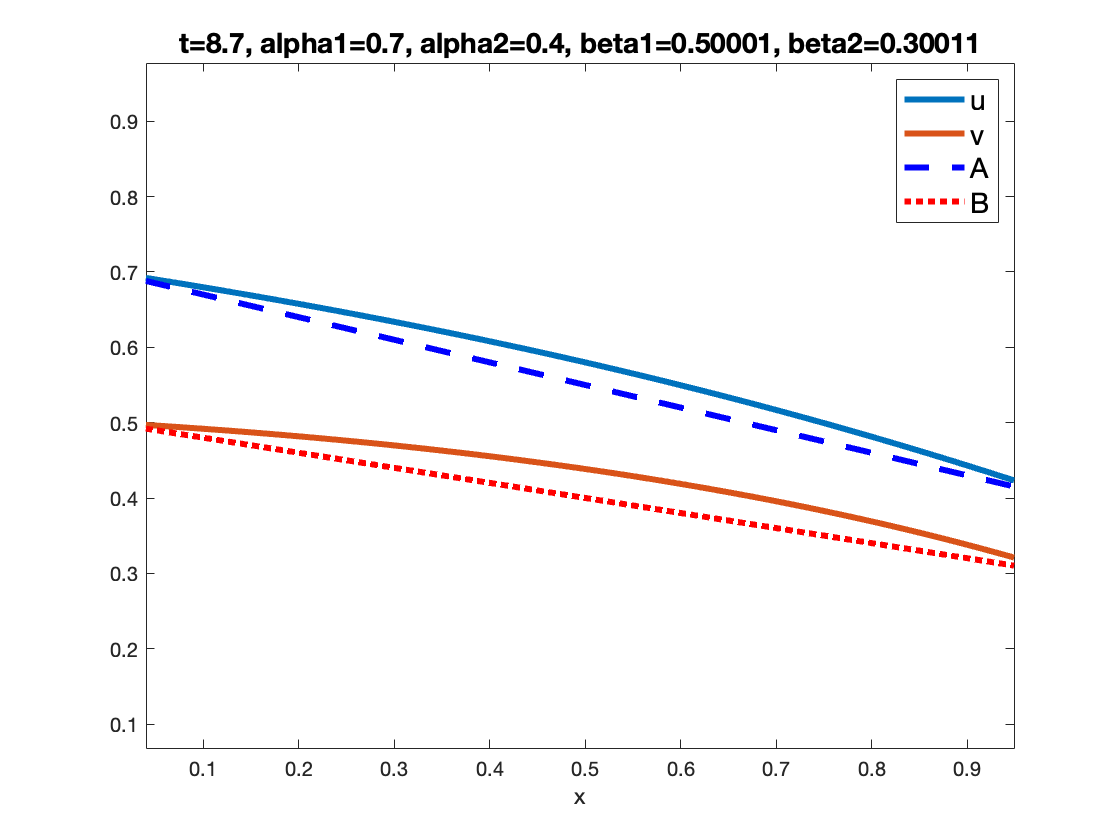}
        \caption{$t=8.7$}
        \label{fig5}
    \end{subfigure}
    \qquad
    \begin{subfigure}{0.45\textwidth}
        \includegraphics[width=\textwidth]{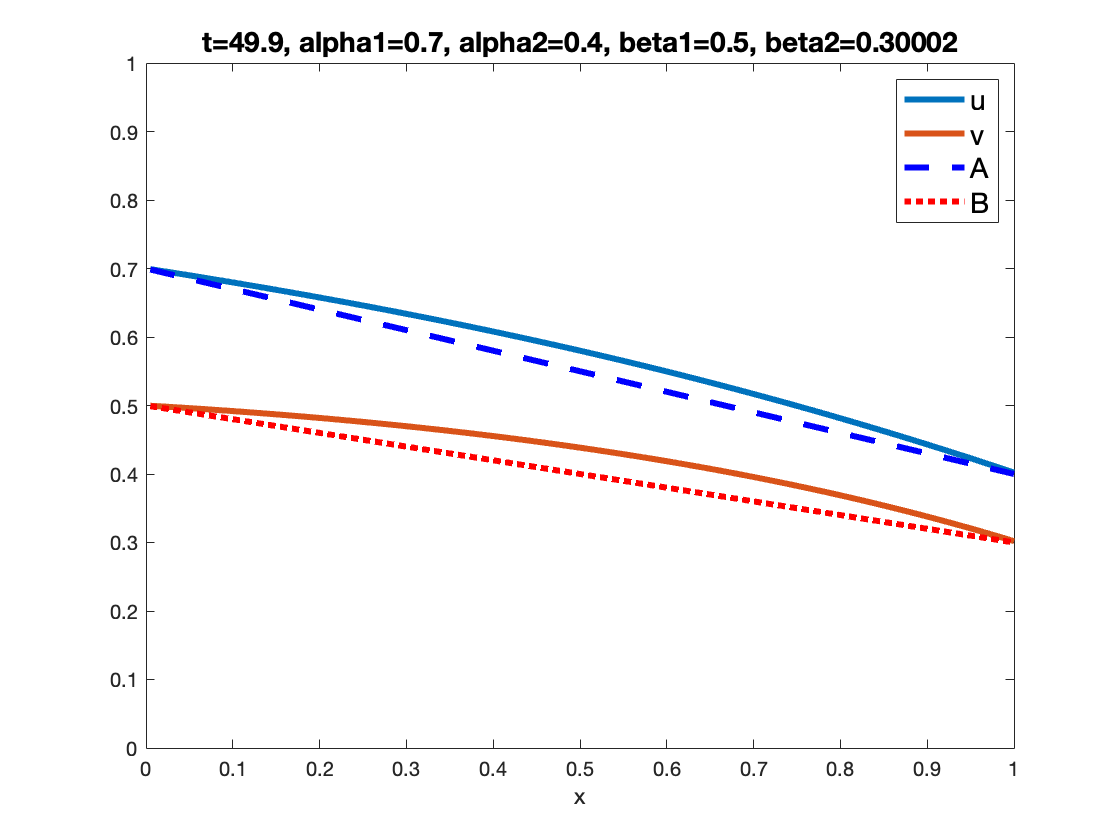}
        \caption{$t=49.9$}
        \label{fig6}
    \end{subfigure}
    \caption{Solution of \eqref{1.1} with $\alpha_1\to\bar\alpha_1\neq\bar\alpha_2\gets\alpha_2,\ \ \beta_1\to\bar\beta_1\neq\bar\beta_2\gets\beta_2$.}
    \label{case3}
\end{figure}\\
In this case, we plot $u,v, A, B$ with $\alpha_1=0.7+\exp(-200000t)$, $\alpha_2=0.4+\exp(-200000t)$,  $\beta_1=0.5+1/(1+10000t)$ and $\beta_2= 0.3+1/(4+1000t)$.
Figure \ref{fig5} plots the solution $(u,v)$ and the linear interpolation $(A,B)$ at time $t=8.7$, with  $\alpha_1=0.7, \alpha_2=0.4$, $\beta_1=0.500, \beta_2=0.300$. Figure \ref{fig6} plots $(u,v,A,B)$ at time $t=49.9$, with  $\alpha_1=0.7, \alpha_2=0.4$, $\beta_1=0.500, \beta_2=0.300$. In particular, we observe that the solution reaches a steady state different from $\big((\bar\alpha_2-\bar\alpha_1)x+\bar\alpha_1,(\bar\beta_2-\bar\beta_1)x+\bar\beta_1\big)=\big(-0.3x+0.7,-0.2x+0.5\big)$.

\subsection{Study of solutions  of \eqref{1.2}} This section contains numerical study of steady state solutions of \eqref{1.2} for different conditions on the boundary data. For simulations in this section, we use the boundary data $\alpha_1, \alpha_2$ for $u$, but notice that the boundary conditions for $v$ cannot be imposed. In this case, they are implicit by \eqref{1.2b}. Numerically, this presents a challenge that we have overcome by calculating explicitly the values at the boundary at each step using \eqref{1.2b}. We still make use of the notation for $\beta_1, \beta_2$ for storing these data values. Throughout these calculations $\varepsilon=0$.  

We observe the following. 
\begin{itemize}
\item \textbf{Case 1}  $\boxed{\alpha_1=\alpha_2\to\bar\alpha}$ 
\end{itemize}
The solution $(u, v)$ converges to the steady state $(\bar\alpha, \bar v)$ as proved in  Theorem \ref{thm2}. \\
\begin{figure}[htbp]
    \begin{subfigure}{0.45\textwidth}
        \includegraphics[width=\textwidth]{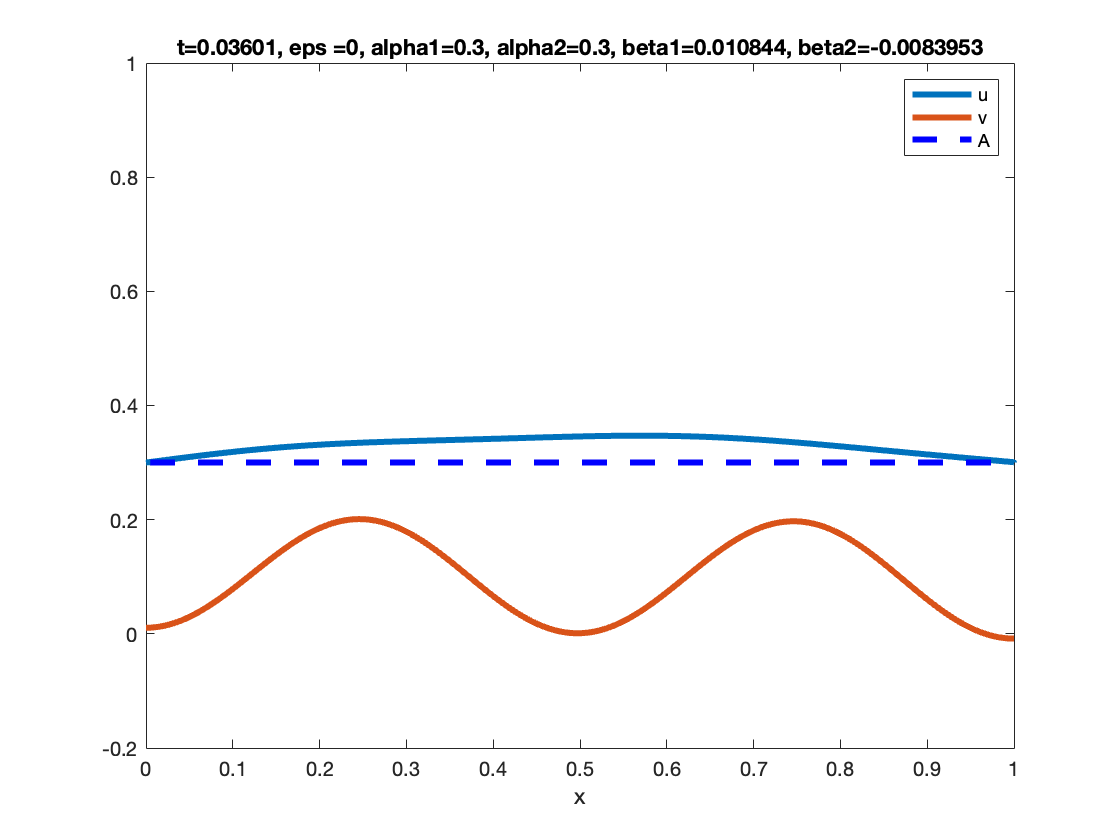}
        \caption{$t=0.036$}
        \label{fig7}
    \end{subfigure}
    \qquad
    \begin{subfigure}{0.45\textwidth}
        \includegraphics[width=\textwidth]{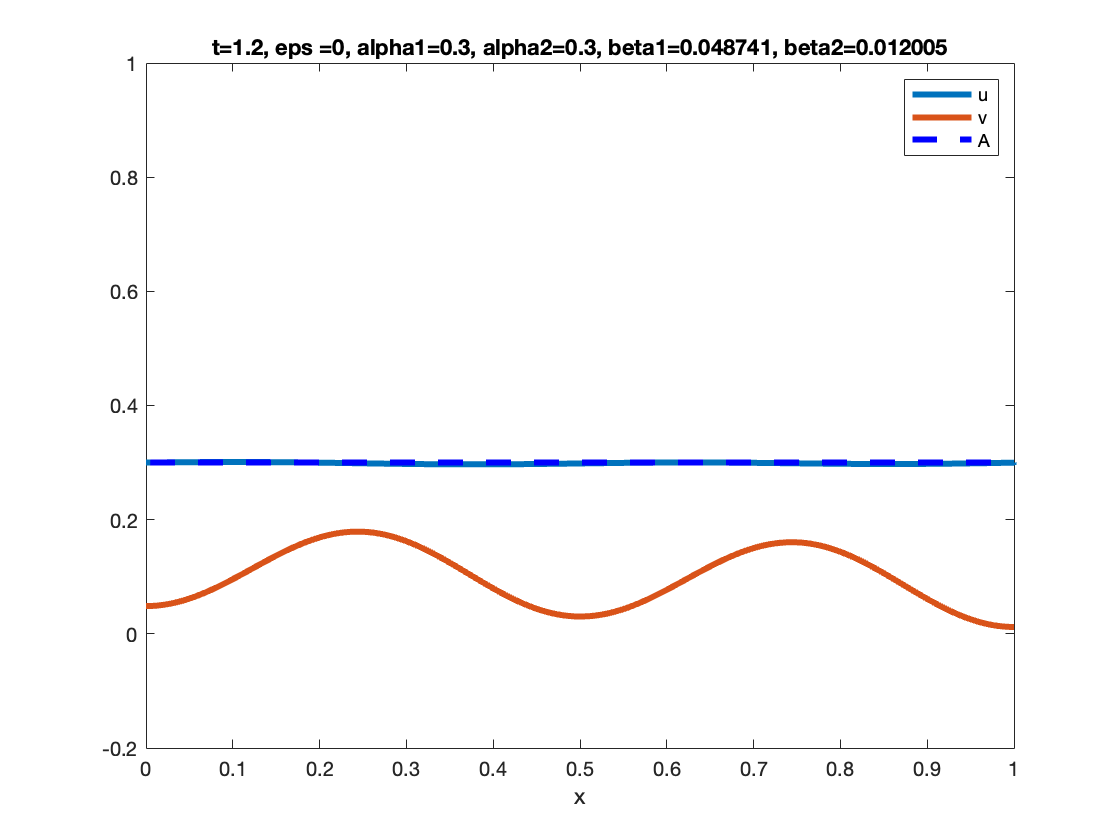}
        \caption{$t=1.2$}
        \label{fig8}
    \end{subfigure}\\
    \begin{subfigure}{0.45\textwidth}
        \includegraphics[width=\textwidth]{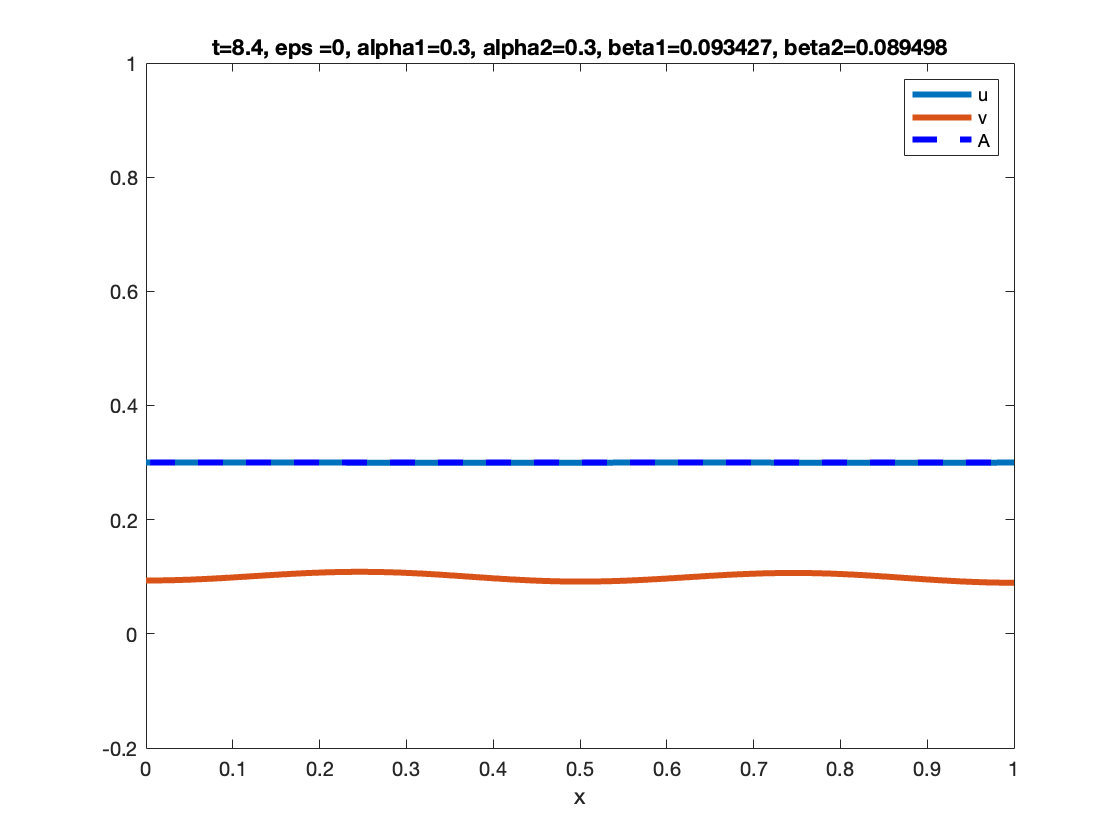}
        \caption{$t=8.4$}
        \label{fig9}
    \end{subfigure}
    \qquad
    \begin{subfigure}{0.45\textwidth}
        \includegraphics[width=\textwidth]{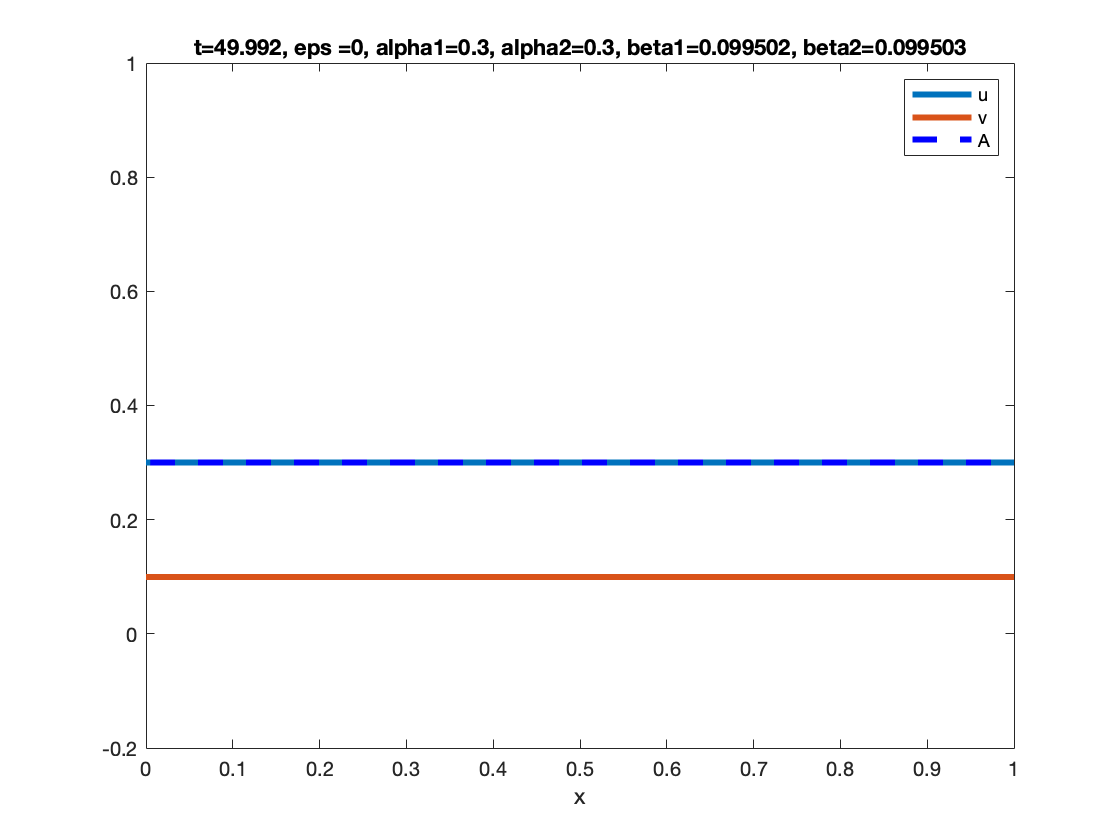}
        \caption{$t=49.99$}
        \label{fig10}
    \end{subfigure}
    \caption{Solution of \eqref{1.2} with $\alpha_1=\alpha_2\to\bar\alpha$.}
     \label{case4}
\end{figure}\\
In this case, we plot $u,v, A$ with $\alpha_1=0.3+\exp(-200000t)=\alpha_2$. Notice that the average of $v_0$, $\bar v=0.1$. Figure \ref{case4} shows the evolution in time of the solution $(u,v)$ and the linear interpolation $A$ at times $t=0.036$, $t=1.2$, $t=8.4$ and $t=49.99$. In particular, we observe that the solution approaches the steady state $(A,\bar v)=(0.3,0.1)$, as $t\to\infty$. Notice that even if the chemical diffusion coefficient is zero, the function $v$ still smooths up and converges to its initial average, as predicted in Theorem \ref{thm2}.

\begin{itemize}
\item \textbf{Case 2}  $\boxed{\alpha_1\neq \alpha_2,\ \  \alpha_1\to\bar\alpha\gets\alpha_2}$ 
\end{itemize}
\padi{In this case, $u$ converges to $\bar\alpha$ and $v$ converges to a steady--state different from $\bar v$.}\\
\begin{figure}[htbp]
    \begin{subfigure}{0.45\textwidth}
        \includegraphics[width=\textwidth]{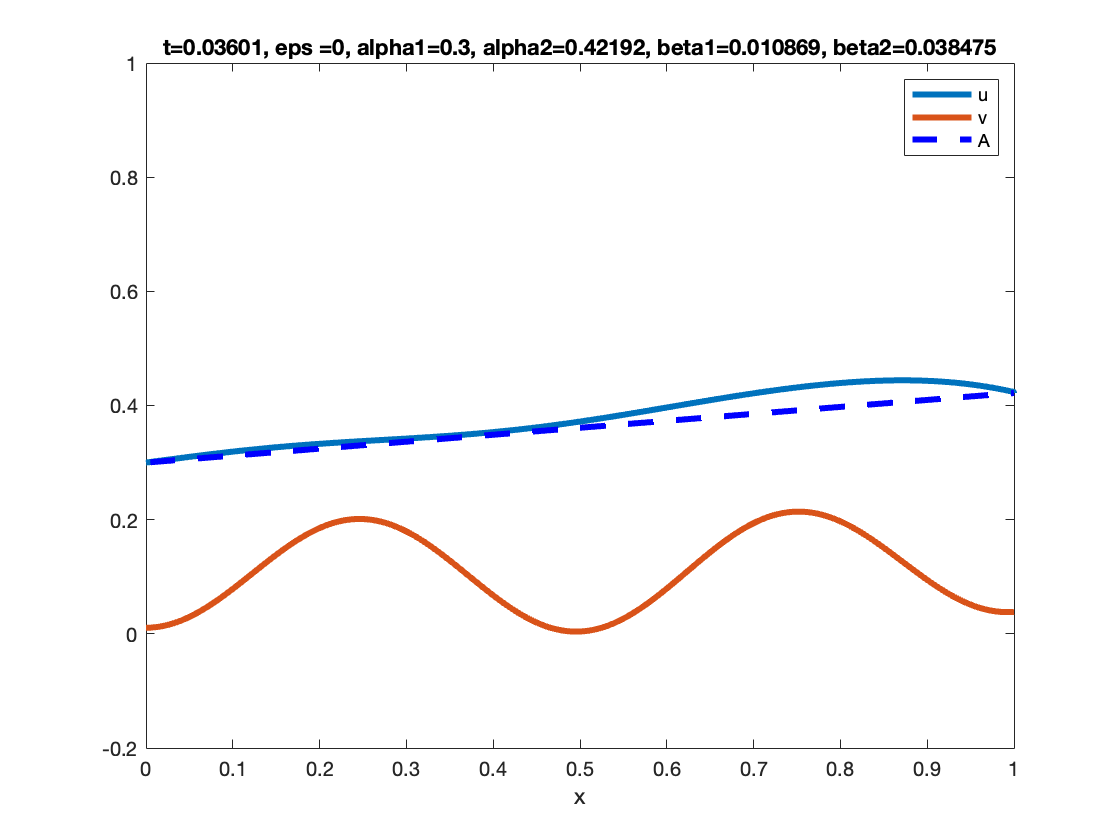}
        \caption{$t=0.036$}
        \label{fig15}
    \end{subfigure}
    \qquad
    \begin{subfigure}{0.45\textwidth}
        \includegraphics[width=\textwidth]{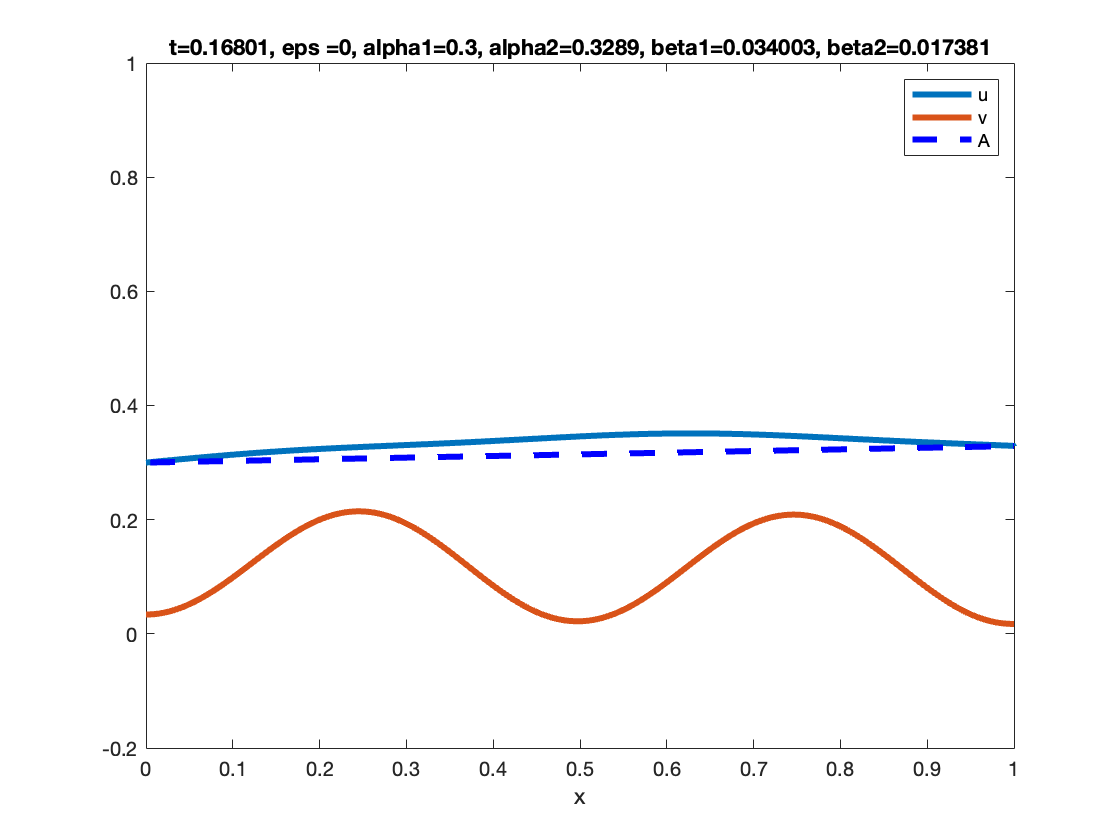}
        \caption{$t=0.168$}
        \label{fig16}
    \end{subfigure}\\
    \begin{subfigure}{0.45\textwidth}
        \includegraphics[width=\textwidth]{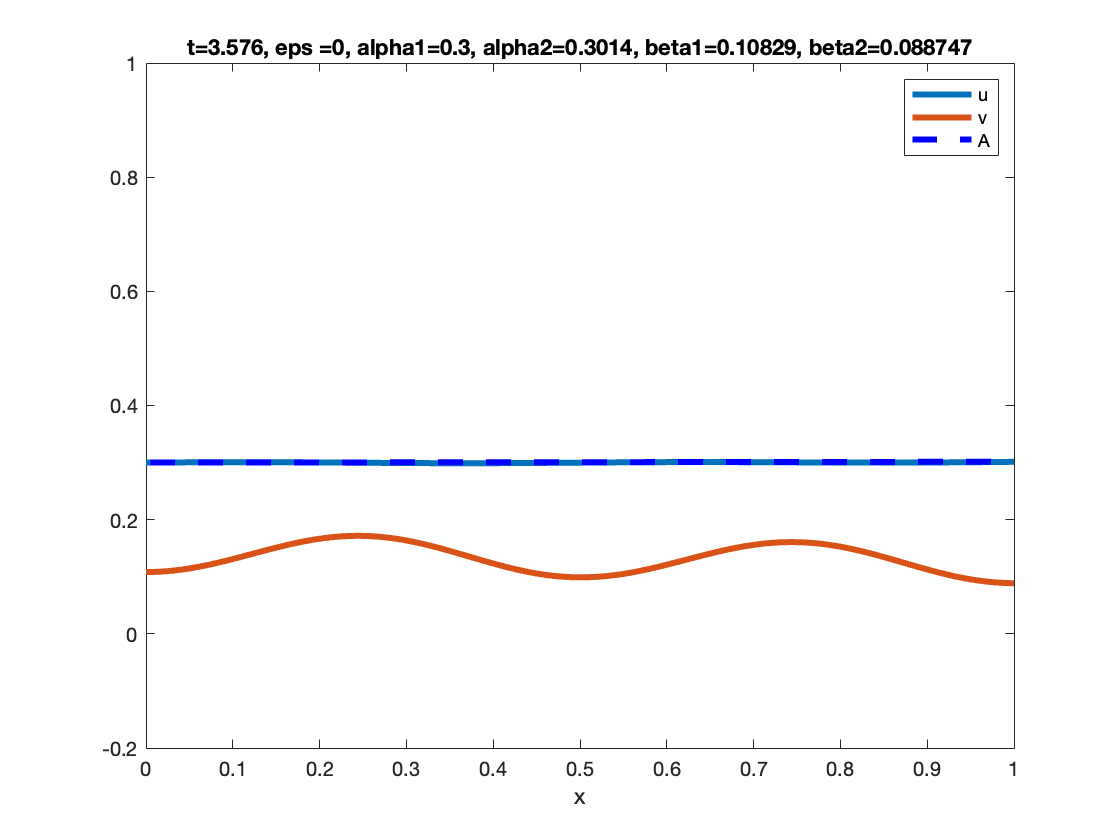}
        \caption{$t=3.576$}
        \label{fig17}
    \end{subfigure}
    \qquad
    \begin{subfigure}{0.45\textwidth}
        \includegraphics[width=\textwidth]{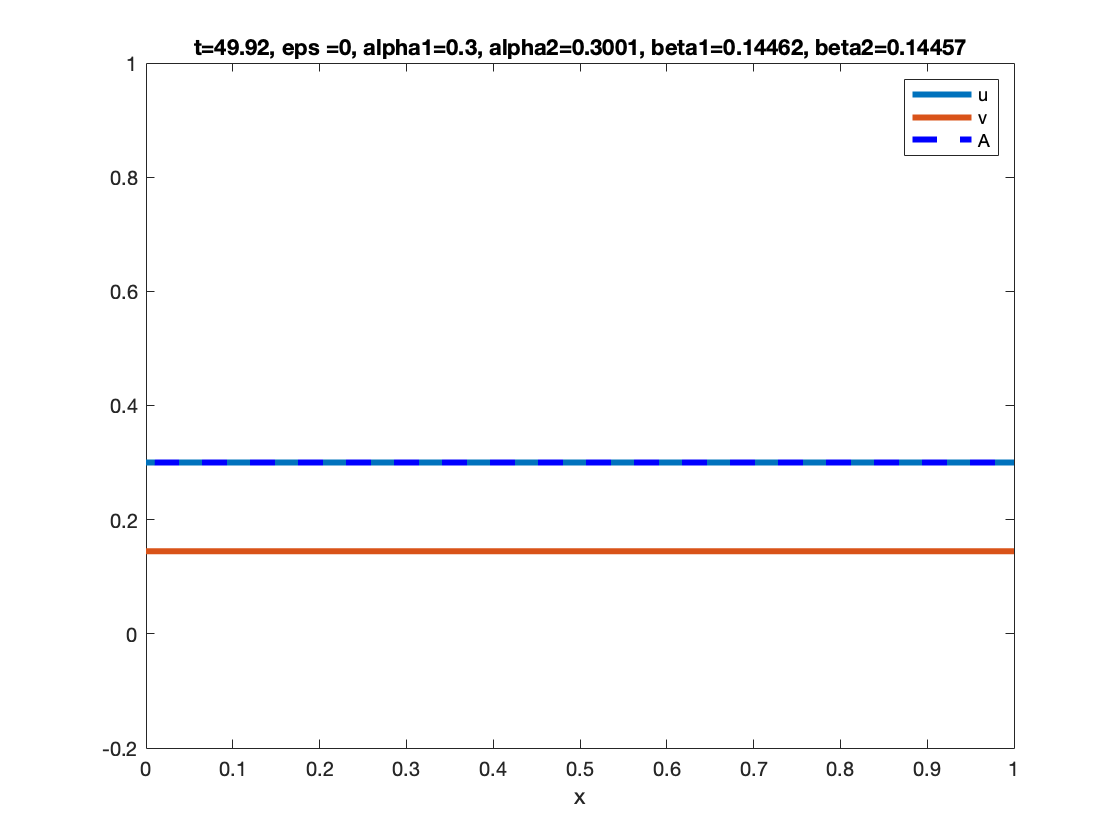}
        \caption{$t=49.92$}
        \label{fig18}
    \end{subfigure}
    \caption{Solution of \eqref{1.2} with $\alpha_1\neq\alpha_2$, \ \   $\alpha_1\to\bar\alpha\gets\alpha_2$.}
    \label{case6}
\end{figure}\\
We plot $u,v, A$ with $\alpha_1=0.3-\exp(-200000t)$ and $\alpha_2=0.3+1/(1+200t)$. Figure \ref{case6} shows the evolution in time of the solution $(u,v)$ at times $t=0.036$, $t=0.168$, $t=3.576$, and $t=49.92$. \padi{In particular, we observe that $u$ converges to the steady state $\bar\alpha=0.3$, but $v$ converges to a steady--state different than $\bar v=0.1$. This shows that the assumptions in Theorem \ref{thm2} are in fact necessary.}

\begin{itemize}
\item \textbf{Case 3}  $\boxed{\alpha_1\to\bar\alpha_1\neq\bar\alpha_2\gets\alpha_2}$ 
\end{itemize}
The solution $u$ converges to a steady state different from $(\bar\alpha_2-\bar\alpha_1)x+\bar\alpha_1$, but $v$ keeps growing or decaying depending on $\alpha_2-\alpha_1$ being negative or positive respectively. \\
\begin{figure}[htbp]
    \begin{subfigure}{0.45\textwidth}
        \includegraphics[width=\textwidth]{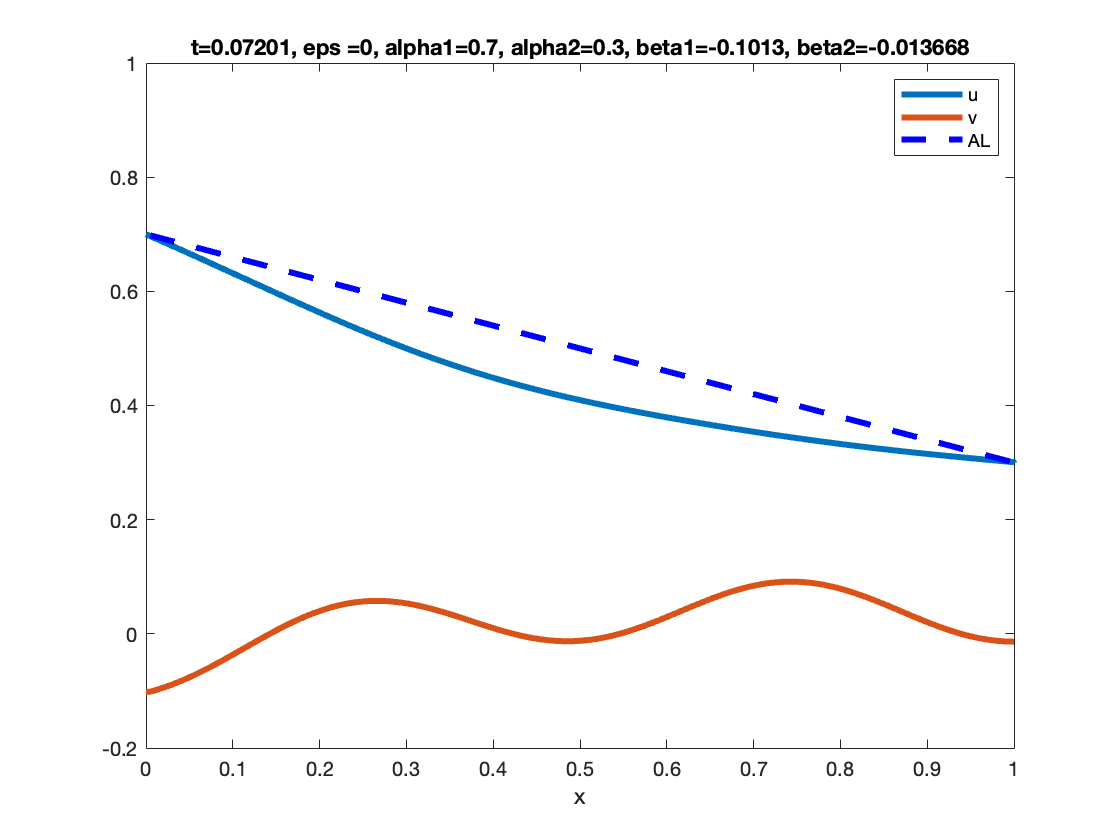}
        \caption{$t=0.072$}
        \label{fig11}
    \end{subfigure}
    \qquad
    \begin{subfigure}{0.45\textwidth}
        \includegraphics[width=\textwidth]{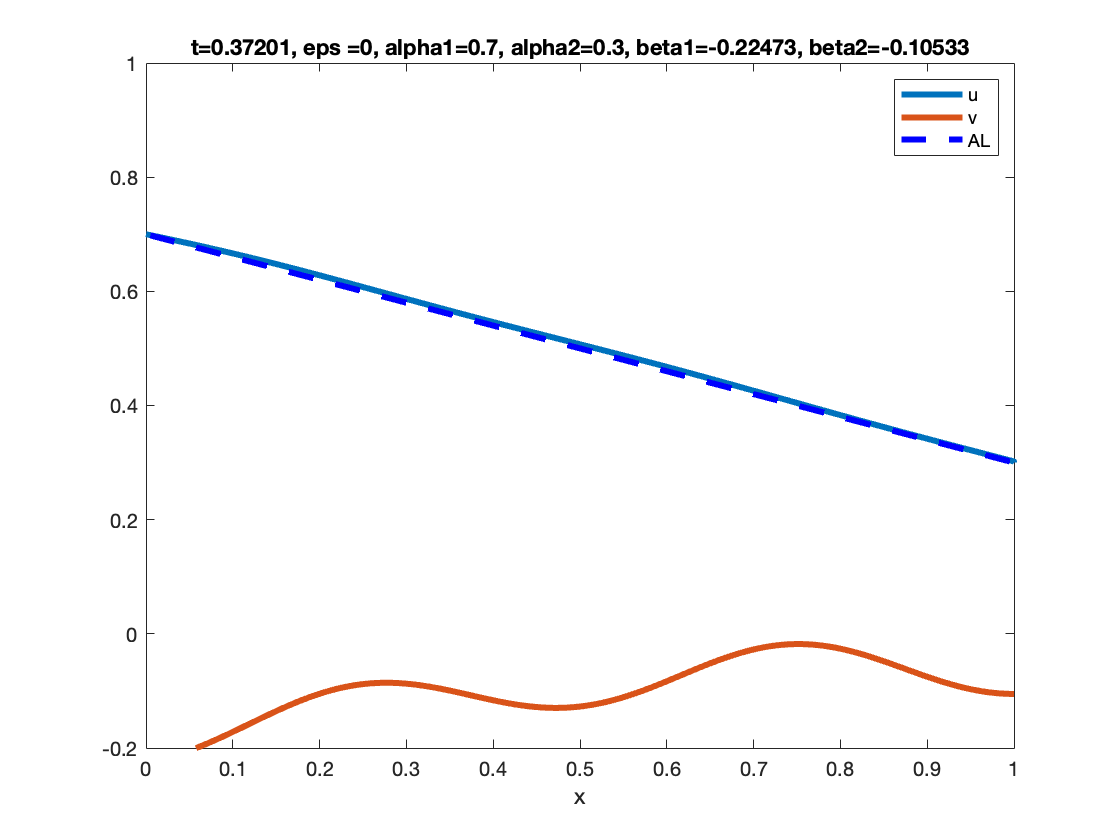}
        \caption{$t=0.372$}
        \label{fig12}
    \end{subfigure}\\
    \begin{subfigure}{0.45\textwidth}
        \includegraphics[width=\textwidth]{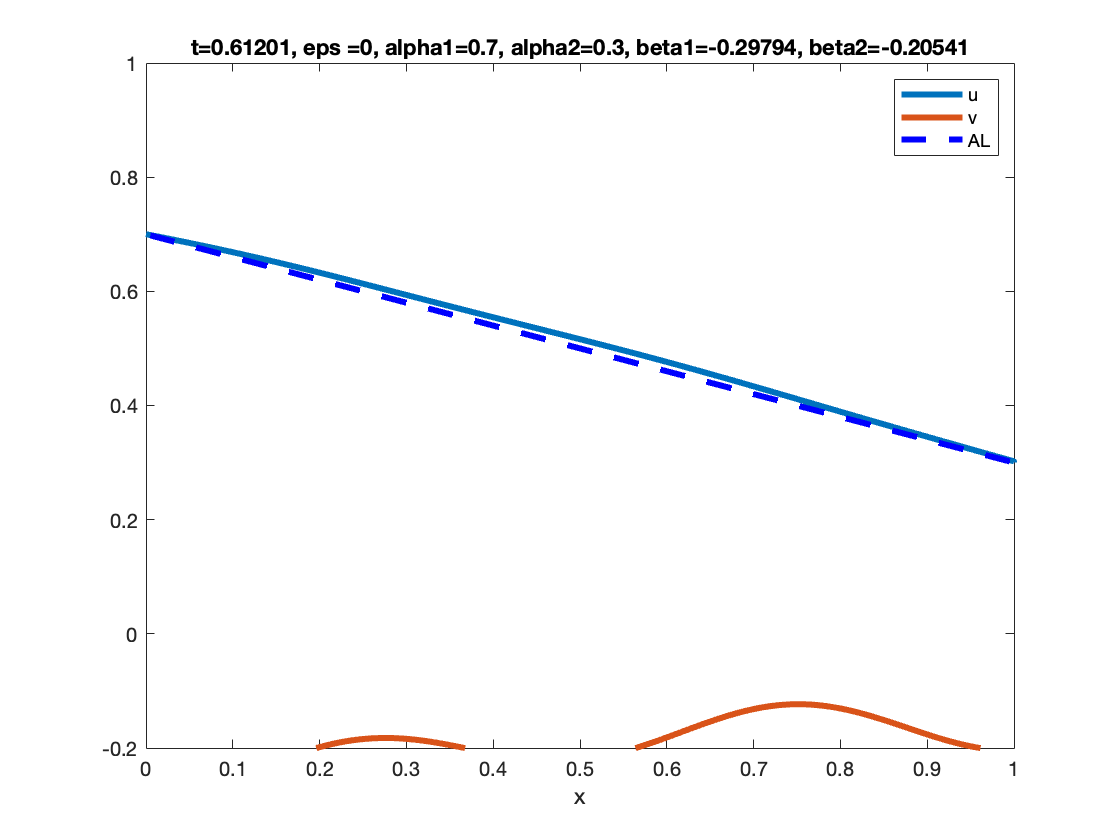}
        \caption{$t=0.612$}
        \label{fig13}
    \end{subfigure}
    \qquad
    \begin{subfigure}{0.45\textwidth}
        \includegraphics[width=\textwidth]{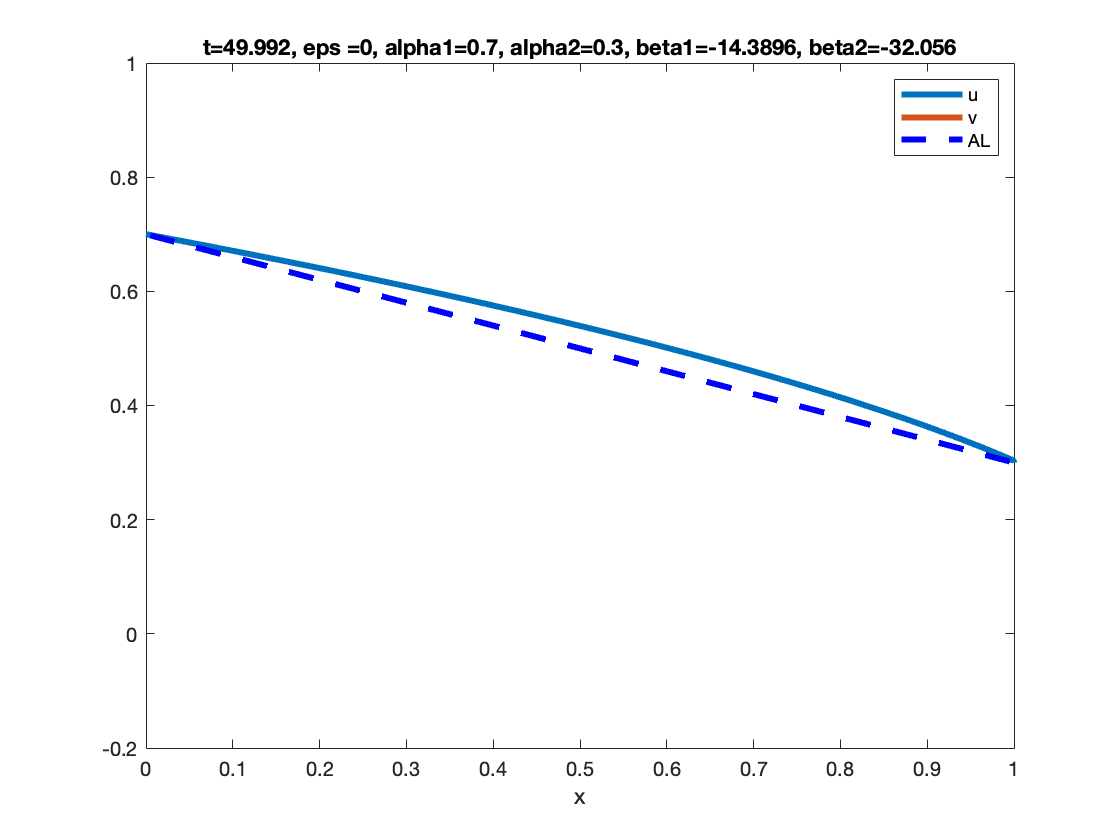}
        \caption{$t=49.99$}
        \label{fig14}
    \end{subfigure}
    \caption{Solution of \eqref{1.2} with $\alpha_1\to\bar\alpha_1\neq\bar\alpha_2\gets\alpha_2$.}
    \label{case5}
\end{figure}\\
In this case, we plot $u,v, A$ with $\alpha_1=0.7+\exp(-200000t)$ and $\alpha_2=0.3+\exp(-200000t)$. Figure \ref{case5} shows the evolution in time of the solution $(u,v)$ at times $t=0.072$, $t=0.372$, $t=0.612$ and $t=49.99$. In particular, we observe that $u$ reaches a steady state different $(\bar\alpha_2-\bar\alpha_1)x+\bar\alpha_1=-0.4x+0.7$, and that $v$ diverges driven by \eqref{1.2b}. 
\section{Conclusion}

We have studied the dynamical behavior of classical solutions to the system of balance laws \eqref{TOS1} with $\chi\mu>0$, derived from the Keller-Segel type model of chemotaxis with logarithmic sensitivity \eqref{OS}, on a finite interval subject to time-dependent Dirichlet type boundary conditions. For the model with $\varepsilon>0$, it is shown that under the assumptions \eqref{aa} and \eqref{ab}, classical solutions are globally well-posed and the perturbations around the reference profiles interpolating the boundary data converge to zero as time goes to infinity. There is no smallness restriction on the strength of the initial perturbations (see Theorem \ref{thm1}). This result generalized previous ones in the sense that the values of the function $v$ at the endpoints of the spatial interval are allowed to be different at any time. When $\varepsilon=0$, similar results are obtained, except the steady state of $v$ is given by its initial average over the spatial interval (see Theorem \ref{thm2}). 

Numerically, we have seen that by keeping the asymptotic convergence of the boundary data while relaxing the assumptions of Theorems \ref{thm1} and \ref{thm2}, the solutions still converge to steady states, except in the case of zero chemical diffusion and unequal end-states of the boundary data for $u$, in which $v$ diverges. For convergent  solutions, when the end-states of the boundary data of either $u$ and/or $v$ are different, the steady states are non-trivial functions of $x$ and are different from the linear profiles interpolating the corresponding end-states of the boundary data. 

At the same time, several of the numerical results observed in this paper still remain to be proved analytically. For example, whether the matched boundary condition for the function $u$ can be relaxed so that similar results as Theorems \ref{thm1}-\ref{thm2} can be obtained (see Remark \ref{rem3}). Another direction of future exploration is the dynamical behavior of the appended model with logistic growth in the $u$-equation. Because of the enhanced dissipation mechanism induced by logistic damping, the chemotaxis-growth model may give us a definite answer to the preceding question. 

\vspace{.2 in}

\noindent{\bf Acknowledgements.} The research of K. Zhao was partially supported by the Simons Foundation's Collaboration Grant for Mathematicians No. 413028. \padi{The authors would like to thank Dr. Ricardo Cortez for the guidance in the numerical simulations.}

\end{document}